\documentclass[11pt]{article}
\usepackage{fullpage}
\usepackage{graphicx}
\usepackage{amsmath,amssymb}
\usepackage{epstopdf}

\usepackage{fullpage}
\usepackage{color}

\renewcommand{\mod}{ {\,{\rm  mod}\;}}

\newcommand{\Ore}{DHGO}
\def\qed{\hfill
\ifhmode\unskip\nobreak\fi\quad\ifmmode\Box\else$\Box$\fi\\ }
\newtheorem{theorem}{Theorem}
\newtheorem{cor}[theorem]{Corollary}
\newtheorem{lemma}[theorem]{Lemma}
\newtheorem{defn}[theorem]{Definition}

\newtheorem{claim}[theorem]{Claim}
\newtheorem{fact}[theorem]{Fact}
\newtheorem{conj}[theorem]{Conjecture}
\newtheorem{const}[theorem]{Construction}
\title{A Brooks-type result for sparse critical graphs}

\author{Alexandr  Kostochka\thanks{
University of Illinois at Urbana--Champaign,  Urbana, IL 61801, USA and
 Sobolev Institute of Mathematics,  Novosibirsk 630090, Russia. Email:
 kostochk@math.uiuc.edu.
Research of this author
is supported in part by NSF grants  DMS-1266016 and  DMS-1600592 and by
grants 15-01-05867 and 16-01-00499  of the Russian Foundation for Basic Research. }
\and
Matthew Yancey\thanks{Department of Mathematics, University of Illinois, Urbana,
IL 61801, USA. E-mail: yancey1@illinois.edu.
Research of this author is partially supported by the
Arnold O. Beckman Research Award of the University of Illinois
at Urbana-Champaign.}}

\begin{document}
 \maketitle

\begin{abstract}
A graph $G$ is  $k$-{\em critical} if it has chromatic number $k$, but every
proper subgraph of $G$ is $(k-1)$--colorable.
Let $f_k(n)$ denote the minimum number of edges in an $n$-vertex $k$-critical graph.
Recently the authors
 gave a lower bound, $f_k(n)  \geq \left\lceil \frac{(k+1)(k-2)|V(G)|-k(k-3)}{2(k-1)}\right\rceil$,
 that solves a conjecture by Gallai from 1963 and is sharp for every $n\equiv 1\,({\rm mod }\, k-1)$.
It is also sharp for $k=4$ and every $n\geq 6$.
In this paper we refine the result by describing
all $n$-vertex $k$-critical graphs $G$ with $|E(G)|= \frac{(k+1)(k-2)|V(G)|-k(k-3)}{2(k-1)}$.
In particular, this result implies exact values of $f_5(n)$ for $n\geq 7$.
\end{abstract}

\section{Introduction}

A {\em proper $k$-coloring}, or simply $k$-{\em coloring}, of a graph $G = (V, E)$ is a function $f:V \rightarrow \{1,2,\dots,k\}$
such that for each $uv \in E$, $f(u) \neq f(v)$.
A graph $G$ is $k$-{\em colorable} if there exists a $k$-coloring of $G$.
The {\em chromatic number}, $\chi(G)$, of a graph $G$ is the smallest $k$ such that $G$ is $k$-colorable.
A graph $G$ is $k$-{\em chromatic} if $\chi(G)=k$.

A graph $G$ is $k$-{\em critical} if $G$ is $k$-chromatic, but every proper subgraph of $G$ is $(k-1)$-colorable.
Critical graphs were first defined and used by Dirac~{\cite{D0,D02,D03}} in 1951-52.
A reason to study $k$-critical graphs is that every $k$-chromatic graph contains a $k$-critical subgraph and
$k$-critical graphs have more restricted structure.
For example, $k$-critical graphs are $2$-connected and $(k-1)$-edge-connected.

One of the basic questions on $k$-critical graphs is: {\em What is the minimum number
$f_k(n)$ of edges in a $k$-critical graph with $n$ vertices?}
This question was first asked by Dirac~\cite{D1} in 1957 and then was reiterated by Gallai~\cite{G2} in 1963, Ore~\cite{O}
in 1967 and
others~\cite{J, J2,Tuza1}.
Gallai~\cite{G2} has found the values of $f_k(n)$ for $n\leq 2k-1$.

\begin{theorem}[Gallai~\cite{G2}] \label{gallai1}
If $k \geq 4$ and $k+2\leq n\leq 2k-1$, then
$$f_k(n)=\frac{1}{2}  \left((k-1)n+(n-k)(2k-n)\right)-1.$$
 \end{theorem}

Kostochka and Stiebitz~\cite{K2} found the value $f_k(2k) = k^2 - 3$.
Gallai~\cite{G1} also conjectured the exact value for $f_k(n)$ for $n\equiv 1\,(\mod k-1)$.

\begin{conj}[Gallai~\cite{G1}] \label{gallai conj}
If $k \geq 4$ and $n\equiv 1 \,(\mod k-1)$, then
$$f_k(n) = \frac{(k+1)(k-2)n-k(k-3)}{2(k-1)} .$$
 \end{conj}

The upper bound on $f_k(n)$ follows from Gallai's construction of $k$-critical graphs with only one vertex of degree at least $k$.
So the main difficulty of the conjecture is in proving the lower bound on $f_k$.

For a graph $G$ and a vertex $u \in V(G)$, a \emph{split} of $u$ is a construction of a new graph $G'$ such that $V(G') = V(G) - u +\{u',u''\}$, where $G - u \cong G' - \{u',u''\}$, $N(u') \cup N(u'') = N(u)$, and $N(u') \cap N(u'') = \emptyset$.
A \emph{\Ore-composition} $O(G_1,G_2)$ of graphs $G_1$ and $G_2$ is a graph obtained as follows:
delete some edge $xy$ from $G_1$,
split some vertex $z$ of $G_2$ into two vertices $z_1$ and $z_2$ of positive degree, and
identify $x$ with $z_1$ and $y$ with $z_2$.  
Note that \Ore-composition could be found in  Dirac's paper~\cite{Dirac64}
and has roots in \cite{Dirac53}. It was also used by Gallai~\cite{G1}
and Haj\'{o}s~\cite{Hajos61}. Ore \cite{O} used it for a composition of complete
graphs.

The mentioned authors observed that if $G_1$ and $G_2$ are $k$-critical and $G_2$ is not $k$-critical after $z$ has been split, then $O(G_1,G_2)$ also is $k$-critical.
This observation implies
\begin{equation}\label{upper f_k}
f_k(n + k - 1) \leq f_k(n) + \frac{(k+1)(k-2)}2 = f_k(n) + (k-1)\frac{(k+1)(k-2)}{2(k-1)}.
\end{equation}
Ore believed that  this construction starting from an extremal graph on at most $2k$ vertices with $G_2 = K_k$ at each iteration yields  sparsest  $k$-critical $n$-vertex graphs for all $n\geq 2k$.
\begin{conj} [Ore~\cite{O}] \label{Ore Conj}
If $k \geq 4$, $n\geq k$ and $n\neq k+1$, then
$ f_k(n + k - 1) = f_k(n) + (k-2)(k+1)/2. $
\end{conj}
Note that Conjecture \ref{gallai conj} is equivalent to the case  $n\equiv 1 \,(\mod k-1)$ of Conjecture \ref{Ore Conj}.

Some lower bounds on $f_k(n)$ were obtained in~\cite{D1,Kr2,G1,K2,K5,FM}.
Recently, the authors~\cite{KY} proved Conjecture \ref{gallai conj} valid.
\begin{theorem} [\cite{KY}] \label{k-critical}
If $k \geq 4$ and $G$ is $k$-critical, then $ |E(G)| \geq  \left\lceil \frac{(k+1)(k-2)|V(G)|-k(k-3)}{2(k-1)}\right\rceil$.
In other words, if $k\geq 4$ and $n\geq k,\,n\neq k+1$, then
$$f_k(n)\geq  F(k,n):=\left\lceil \frac{(k+1)(k-2)n-k(k-3)}{2(k-1)}\right\rceil.$$
\end{theorem}

The result also confirms Conjecture \ref{Ore Conj} in several cases.

\begin{cor} [\cite{KY}] \label{Ore Cor}
Conjecture \ref{Ore Conj} is true if
(i) $k=4$,
	(ii) $k=5$ and $n \equiv 2 \,(\mod 4)$, or
	(iii) $n \equiv 1 \,(\mod k-1)$.
\end{cor}

Some applications of Theorem~\ref{k-critical} are given in~\cite{KY} and~\cite{BKLY}. In~\cite{KY2},
the authors derive from a partial case of Theorem~\ref{k-critical} a half-page proof of the
well-known Gr\" otzsch Theorem~\cite{Gr} that every planar triangle-free graph is 3-colorable.
Conjecture \ref{Ore Conj} is still open in general.
By examining known values of $f_k(n)$ when $n \leq 2k$, it follows that $f_k(n) - F(k,n) \leq k^2/8$.

The goal of this paper is to describe the {\em $k$-extremal} graphs, i.e. the
$k$-critical graphs $G$  such that $|E(G)|=\frac{(k+1)(k-2)|V(G)|-k(k-3)}{2(k-1)}$. This is
a refinement of Conjecture \ref{gallai conj}: For $n\equiv 1 \,(\mod k-1)$, we describe all $n$-vertex
$k$-critical graphs $G$ with $|E(G)|=f_k(n)$.
This is also the next step towards  the full solution of Conjecture \ref{Ore Conj}.

By definition, if $G$ is $k$-extremal, then $\frac{(k+1)(k-2)|V(G)|-k(k-3)}{2(k-1)}$ is an integer, and so $|V(G)|\equiv 1 \,(\mod k-1)$.
For example, $K_{k}$ is $k$-extremal.

Suppose that $G_1$ and $G_2$ are $k$-extremal and $G=O(G_1,G_2)$.
Then
$$|E(G)|=|E(G_1)|+|E(G_2)|-1=\frac{(k+1)(k-2)(|V(G_1)|+|V(G_2)|)-2k(k-3)}{2(k-1)}-1
$$
$$=\frac{(k+1)(k-2)|V(G)|-k(k-3)}{2(k-1)}.$$
After $z$ is split, $G_2$ will still have $F(k, |V(G_2)|) < F(k, |V(G_2)|+1)$ edges, and therefore will not be $k$-critical.
Thus the \Ore-composition of any two $k$-extremal graphs is again $k$-extremal.

A graph is a $k$-{\em Ore graph} if it is obtained from a set of copies of $K_k$ by a sequence of \Ore-compositions.
By the above, every $k$-Ore graph is $k$-extremal.
This yields an explicit construction of infinitely many $k$-extremal graphs.

The main result of the present paper is the following.

\begin{theorem} \label{ext}
Let $k \geq 4$ and $G$ be a $k$-critical graph. Then $G$ is $k$-extremal if and only if it
 is a $k$-Ore graph. Moreover, if $G$ is not a $k$-Ore graph, then
  $|E(G)|\geq\frac{(k+1)(k-2)|V(G)|-y_k}{2(k-1)}$,
 where $y_k=\max\{2k-6,k^2-5k+2\}$.
Thus $y_4 = 2$, $y_5 = 4$, and $y_k = k^2 - 5k + 2$ for $k \geq 6$.
\end{theorem}

The message of Theorem~\ref{ext} is that although for every $k\geq 4$ there are infinitely many $k$-extremal graphs, they all have a simple structure.
In particular, every $k$-extremal graph distinct from $K_k$ has a separating set of size $2$.
The theorem gives a slightly better approximation for $f_k(n)$ and adds new cases for which we now know the exact values of $f_k(n)$:

\begin{cor} \label{new tightness}
Conjecture \ref{Ore Conj} holds and the value of  $f_k(n)$ is known if
(i)  $k\in\{4,5\}$,\;
	(ii) $k=6$ and $n \equiv 0 \,(\mod 5)$,\;
	(iii) $k=6$ and $n \equiv 2 \,(\mod 5)$,\;
	(iv) $k=7$ and $n \equiv 2 \,(\mod 6)$, or\;
	(v) $k\geq 4$ and $n \equiv 1 \,(\mod k-1)$.
\end{cor}

This value of  $y_k$ in Theorem~\ref{ext} is best possible in the  sense
 that for every $k\geq 4$,
there exist  infinitely many $3$-connected graphs $G$ with $|E(G)|=\frac{(k+1)(k-2)|V(G)|-y_k}{2(k-1)}$.
The idea of this construction (Construction~\ref{3ConnConst}) and the examples for $k=4,5$ are due to
Toft~(\cite{T12}, based on~\cite{Toft2}). Construction~\ref{con2} produces the examples for $k\geq 6$.

 Theorem~\ref{ext} has already found interesting applications. In~\cite{BDKLY}, it was used to describe the $4$-critical planar graphs with
 exactly $4$ triangles. This problem was studied by Axenov~\cite{aks76} in the seventies, and then mentioned by 
 Steinberg~\cite{steinberg93} (quoting Erd\H os from 1990), and Borodin~\cite{borodinsurvey}. It was proved in~\cite{BDKLY} that
 the $4$-critical planar graphs with
 exactly $4$ triangles and no $4$-faces are exactly the $4$-Ore graphs with exactly $4$ triangles. Also, Kierstead and Rabern~\cite{KR} and
 independently Postle~\cite{Po} have used Theorem~\ref{ext} to describe the infinite family of $4$-critical graphs $G$ with the property that
 for each edge $xy\in E(G)$, $d(x)+d(y)\leq 7$. It turned out that such graphs form a subfamily of the family of $4$-Ore graphs.
 
Our proofs will use the language of \emph{potentials}.

\begin{defn} Let  $G$ be a graph.
For $R \subseteq V(G)$, define \emph{the $k$-potential of $R$} to be
\begin{equation}\label{rho}
\rho_{k,G}(R) =(k+1)(k-2)|R| - 2(k-1)|E(G[R])|.\end{equation}
When there is no chance for confusion, we will use $\rho_k(R)$.
Let $P_k(G) = \min_{\emptyset \neq R \subseteq V(G)} \rho_k(R)$.
\end{defn}

Informally, $\rho_{k,G}(R)$ measures  how many edges are needed to be added to $G[R]$ (or removed, if the potential is negative)
 in order to obtain a graph  with average degree $\frac{(k+1)(k-2)}{k-1}$.
Our proofs below will involve adding and deleting edges and vertices, so using the language of potentials
 helps keep track of whether or not the resulting graph maintains the assumptions of the theorem.


Translated into the language of potentials, Theorem~\ref{k-critical} sounds as follows.

\begin{cor}[\cite{KY}] \label{k(k-3)}
If $G$ is $k$-critical then $ \rho_{k}(V(G)) \leq k(k-3) $.
In particular, if a graph $H$ satisfies $\rho_{k,H}(S) > k(k-3)$ for all nonempty $S \subseteq V(H)$, then $H$ is $(k-1)$-colorable.
\end{cor}

Similarly, our main result, Theorem~\ref{ext}, is:

\begin{theorem}\label{pot theorem}
If $G$ is $k$-critical and not a $k$-Ore graph, then
 $$\rho_{k}(V(G)) \leq y_k,$$
 where $y_k=\max\{2k-6,k^2-5k+2\}$.
 In particular, if a graph $H$ does not contain a $k$-Ore graph as a subgraph and $ \rho_{k,H}(S) > y_k$ for all nonempty $S \subseteq V(H)$, then $H$ is $(k-1)$-colorable.
 \end{theorem}



Our strategy of the proof (similar to those in~\cite{BKo,BKY,KY,KY2}) is to consider a minimum counter-example $J$ to Theorem~\ref{pot theorem} and 
derive a set of its properties leading to a contradiction. Quite useful  claims will be that all nontrivial proper subsets of $V(J)$
have ``high''  potentials. Important examples of such claims are Claim~\ref{very small} and  Lemma~\ref{small potential} below.    This will help us
to provide $(k-1)$-colorings of subgraphs of $J$ with additional properties. For example, Claim~\ref{very small} will imply
Claim~\ref{very small potential} that adding any edge to a  subgraph $H$ of $J$ with $1<|V(H)|<|V(J)|$ leaves the subgraph $(k-1)$-colorable.
Important new ingredient of the proof is the study in the next section of the properties of $k$-Ore graphs and their colorings.
 In Section 3 we prove basic
properties of our minimum counter-example $J$, including Claim~\ref{very small} mentioned above. Then in Section 4 we introduce and study
properties of {\em clusters} -- sets of vertices of degree $k-1$ in $J$ with the same closed neighborhood. This will allow us to
prove Lemma~\ref{small potential}. Based on this lemma and its corollaries, we prove Theorem~\ref{pot theorem} in Section 5 
using some variations of discharging; the cases of small $k$ will need separate considerations.
In Section 6 we  discuss the sharpness of our result and
in Section 7  --- some algorithmic aspects of it.

\section{Potentials and Ore graphs}

The fact below summarizes useful properties of $\rho_k$ and $y_k$ following directly from the definitions or Corollary \ref{k(k-3)}.

\begin{fact}\label{f1} For the $k$-potential defined by (\ref{rho}), we have
\begin{enumerate}
	\item Potential is submodular: 
\begin{equation}\label{a6}
 \rho_{k}(X \cap Y) + \rho_{k}(X \cup Y) = \rho_{k}(X) + \rho_{k}(Y)-2(k-1)|E_G[X-Y,Y-X]|.
\end{equation}
	\item $\rho_{k}(V(K_1)) = (k+1)(k-2)$.
	\item $\rho_{k}(V(K_2)) = 2(k^2-2k-1)$.	
	\item $\rho_{k}(V(K_{k-1})) = 2(k-2)(k-1)$.
	\item $\rho_{k}(V(K_k)) = k(k-3)$.
	\item If $k \geq 4$, then $\rho_{k}(V(K_k)) \leq \rho_{k}(V(K_1)) \leq \rho_{k}(V(K_{k-1})) \leq \rho_{k}(V(K_2)) \leq \rho_{k}(V(K_i))$ for all $3 \leq i \leq k-2$. 
	\item For any vertex set $S$, $\rho_k(S) \geq \rho_k(K_{|S|})$. In particular, if $1 \leq |S| \leq k-1$, then $\rho_k(S) \geq (k+1)(k-2)$.  If $2 \leq |S| \leq k-1$, then $\rho_{k,G}(S) \geq 2(k-2)(k-1)$ with equality only if $|S|=k-1$ and $G[S]=K_{k-1}$.
	\item $k(k-3) \leq y_k+2k-2< (k+1)(k-2)$.
	\item $\rho_k(A)$ is  even for each $k$ and $A$. 
	\item If $G$ is a graph with a spanning subgraph $H$  that  is a $k$-Ore graph, then $\rho_{k,G}(V(G)) \leq k(k-3)$. If $H = G$, then we have equality.  If $H$ is a proper subgraph of $G$, then $\rho_{k,G}(V(G)) \leq y_k$.
\end{enumerate}
\end{fact}


A common technique in constructing critical graphs (see \cite{J,steinberg93}) is to use \emph{quasi-edges} and \emph{quasi-vertices}.
For $k\geq 3$, a graph $G$, and $x,y\in V(G)$, a $k$-{\em quasi-$xy$-edge} $Q_k(x,y)$ is a subset $Q$ of $V(G)$ such that $x,y\in Q$ and\\
 (Q1)  $G[Q]$ has a  $(k-1)$-coloring,\\
 (Q2) $\phi(x)\neq\phi(y)$ for every proper $(k-1)$-coloring of $G[Q]$, and\\
 (Q3) for any edge $e \in G[Q]$, $G[Q] - e$ has a $(k-1)$-coloring $\phi$  such that $\phi(x) = \phi(y)$.\\
Symmetrically, a $k$-{\em quasi-$xy$-vertex} $Q'_k(x,y)$ is a subset $Q'$ of $V(G)$ such that $x,y\in Q'$ and\\
 (Q$'$1)  $G[Q']$ has a  $(k-1)$-coloring,\\
 (Q$'$2) $\phi(x)=\phi(y)$ for every proper $(k-1)$-coloring of $G[Q']$, and\\
 (Q$'$3) for any edge $e \in G[Q']$, $G[Q'] - e$ has a $(k-1)$-coloring $\phi$  such that $\phi(x) \neq \phi(y)$.

If $G$ is a $k$-critical graph, then for each  $e= xy\in E(G)$, graph $G-e$ is a $k$-quasi-$xy$-vertex.
On the other hand, given some $k$-quasi-vertices and $k$-quasi-edges, one can construct from copies of them infinitely many $k$-critical graphs.
In particular, the \Ore-composition can be viewed in this way.
The next observation is well known and almost trivial, but we state it, because we use it often.

 \begin{fact}\label{fa7} Let $k\geq 4$. If a $k$-critical graph $G$ has a separating set $\{x,y\}$, then\\
(1) $G-\{x,y\}$ has exactly two components, say with vertex sets $A'$ and $B'$;\\
(2)  $xy \notin E(G)$;\\
(3) one of  $A'\cup \{x,y\}$ and $B'\cup \{x,y\}$   is a   $k$-quasi-$xy$-edge
and the other is a  $k$-quasi-$xy$-vertex.
\end{fact}

A quasi-edge and a quasi-vertex are very related structures, as seen by the following construction.

\begin{fact}\label{vertex to edge and back}
If $Q_k(x,z)$ is a $k$-quasi-$xz$-vertex and  $Q'(x,y)$ is obtained from $Q_k(x,z)$ by appending a leaf $y$ that is adjacent only to $z$, 
then $Q'(x,y)$ is a $k$-quasi-$xy$-edge.
If $Q_k'(x,y)$ is a $k$-quasi-$xy$-edge and $N(y) \cap Q_k'(x,y) = \{z\}$, then the vertex set $Q_k(x,z) = Q_k'(x,y) - y$ is a $k$-quasi-$xz$-vertex.
\end{fact}

Fact~\ref{fa7} together with the definition of $k$-Ore graphs, implies the following.

 \begin{fact}\label{f2}
Every $k$-Ore graph $G\neq K_k$ has a separating set $\{x,y\}$ and two vertex subsets $A=A(G,x,y)$ and $B=B(G,x,y)$ such that \\
(i) $A\cap B=\{x,y\}$, $A\cup B=V(G)$ and no edge of $G$ connects $A-x-y$ with $B-x-y$,\\
(ii) the graph $ \widetilde G(x,y)$ obtained from $G[A]$ by adding edge $xy$ is a $k$-Ore graph, \\
(iii) the graph $\check G(x,y)$ obtained from $G[B]$ by gluing $x$ with $y$ into a new vertex $x*y$ is a $k$-Ore graph, and \\
(iv) $xy \notin E(G)$.
\end{fact}

In terms of Fact~\ref{f2}, $G$ is the \Ore-composition of $ \widetilde G(x,y)$ and $\check G(x,y)$.
We will repeatedly use the notation in this fact.
The next fact directly follows from the definitions.

\begin{fact}\label{f1.1} Using the notation in Fact \ref{f2}, we have 
\begin{enumerate}
	\item $A$ is a $k${-quasi-$xy$-vertex};
	\item $B$ is a $k${-quasi-$xy$-edge};
	\item $ \rho_{k,G}(A) = \rho_{k,K_1}(V(K_1)) = (k+1)(k-2) $;
	\item $ \rho_{k,G}(B) = \rho_{k,K_2}(V(K_2)) = 2(k^2 - 2k - 1) $;
	\item $N(x) \cap B \cap N(y) = \emptyset$;
	\item $N_{\widetilde G}(v) = N_G(v)$ for each $v \in A - x - y$;
	\item $d_{\check G}(v) = d_G(v)$ for each $v \in B - x - y$;
	\item 
	\begin{itemize}
		\item If $R \subseteq A - x$ or $R \subseteq A - y$, then $\rho_{k,G}(R)=\rho_{k, \widetilde G}(R)$.
		\item If $R \subseteq B - \{x,y\}$, then $\rho_{k,G}(R)=\rho_{k,\check G}(R)$.
		\item If $R \subseteq B - x$ or $R \subseteq B - y$ and $\check R = R - \{x,y\} + x*y$, then $\rho_{k,G}(R)\geq\rho_{k,\check G}(\check R)$.
		\item If $\{x,y\} \subseteq R \subseteq B$ and $\check R = R - \{x,y\} + x*y$, then $\rho_{k,G}(R) = \rho_{k,\check G}(\check R) + (k+1)(k-2)$.
	\end{itemize}
	 \item $\rho_{k,G}(V(G)) = \rho_{k, V(\widetilde G)}(V(\widetilde G)) = \rho_{k, V(\check G)}(V(\check G)) = k(k-3)$.
\end{enumerate}
\end{fact}

\begin{claim} \label{Oresmall} For every $k$-Ore graph $G$ and each $\emptyset\neq R \subsetneq V(G)$, we have $\rho_{k,G}(R) \geq (k+1)(k-2)$.
\end{claim}
{\bf Proof.} Let $G$ be a smallest counter-example to the claim and let
\begin{equation}\label{j14}
 \mbox{
$ R\subsetneq V(G)$ be a smallest  nonempty proper subset of $V(G)$ with 
 $\rho_{k,G}(R)<(k+1)(k-2)$.}
\end{equation}

By Fact~\ref{f1}.7, $G\neq K_k$.
So, by Fact \ref{f2} there is a separating set $\{x,y\}$ and two
vertex subsets $A=A(G,x,y)$ and $B=B(G,x,y)$  as in Fact~\ref{f2}.
By the minimality of $G$, every proper subset of $V( \widetilde  G(x,y))$ and of $V(\check G(x,y))$ has potential at least $(k+1)(k-2)$. 
  If  $G[R]$ were disconnected,
then the vertex set of some component of $G[R]$ would also have potential less than $(k+1)(k-2)$, contradicting the
minimality of $R$.
So, $G[R]$ is connected. Since
$\rho_{k,G}(R)<(k+1)(k-2)$ by~\eqref{j14} and $R$ is non-empty, by Fact~\ref{f1}.7,  $|R| \geq k$.

{\bf Case 1:} $\{x,y\}\cap R=\emptyset$. Since $G[R]$ is connected, $R$ is a
non-empty proper subset either of $A$ or $B$. 
This contradicts Fact \ref{f1.1} and the minimality of $G$.

{\bf Case 2:} $\{x,y\}\cap R=\{x\}$. The set $R\cap A$ induces a non-empty connected subgraph of $G$, and so
by the minimality of $|R|$, $\rho_{k,G}(R\cap A)\geq (k+1)(k-2)$. Similarly,
$\rho_{k,G}(R\cap B)\geq (k+1)(k-2)$. By Fact~\ref{f1}.1,
$$\rho_{k,G}(R)=\rho_{k,G}(R\cap A)+\rho_{k,G}(R\cap B)-\rho_{k,G}(\{x\}) \geq (k+1)(k-2),$$
a contradiction to~\eqref{j14}.

{\bf Case 3:} $\{x,y\}\subseteq R$. If $A\subseteq R$, then by Facts~\ref{f1}.1 and \ref{f1.1}.3,
$$\rho_{k,\check G(x,y)}((R-A)+x*y) = \rho_{k,G}(R) - \rho_{k,G}(A) + \rho_{k, \check G(x,y)}(\{x*y\}) =\rho_{k,G}(R).$$
 But by the minimality of $G$, this
is at least $(k+1)(k-2)$, a contradiction to~\eqref{j14}. Similarly, if $B\subseteq R$, then
$$\rho_{k, \widetilde G(x,y)}(R\cap A)=\rho_{k,G}(R) - \rho_{k,G}(B) + \rho_{k, \widetilde G(x,y)}(\{x,y\}) = \rho_{k,G}(R),$$
 a contradiction to~\eqref{j14} again. So, suppose
$A-R\neq\emptyset$ and $B-R\neq\emptyset$. 
By the minimality of $G$, we have  $\rho_{k, \widetilde G(x,y)}(R\cap A)\geq (k+1)(k-2)$. Since $xy$ is an edge in $ \widetilde G(x,y)$
but not in $G$, this yields $\rho_{k,G}(R\cap A)\geq (k+1)(k-2)+2(k-1)$.
Similarly, $\rho_{k,\check G(x,y)}((R-A)+x*y)\geq (k+1)(k-2)$ and thus
$\rho_{k,G}(R\cap B)\geq 2(k+1)(k-2)$. Then
$$\rho_{k,G}(R)=\rho_{k,G}(R\cap A)+\rho_{k,G}(R\cap B)-2\rho_{k,G}(K_1)\geq (k+1)(k-2)+2(k-1),$$
a contradiction to~\eqref{j14}.
\qed

\begin{defn}\label{def1}
For a set $U$  of vertices in a graph $G$, the  \emph{border }  of $U$ is   $U_*=\{w \in U : N(w) \not\subseteq U\}$, i.e. the
set of vertices in $U$ that have neighbors outside of $U$.
\end{defn}

If $U_*$ is the border  of $U\subset V(G)$ and $U_*\neq U$, then $U_*$ is a separating set in $G$.

\medskip
A set $S$ of vertices in a graph $G$ is {\em standard}, if\\
(a) $\rho_{k,G}(S)=(k+1)(k-2)$ and\\
(b) the border  of $S$ is a $2$-element set  $\{x,y\}$ such that  $G[S-\{x,y\}]$ is connected, and\\
(c) $S$ is a $k$-quasi-$\{x,y\}$-vertex.


Note that a standard set is a  $k$-quasi-vertex whose $k$-potential is the same as that of a vertex.

\begin{lemma}\label{ore4} Let  $G$ be a $k$-Ore graph.
Let $W\subset V(G)$ with $|W|\geq 2$ and $\rho_k(W)\leq(k+1)(k-2)$.
Then $G[W]$ is connected and contains a standard set.
\end{lemma}

{\bf Proof.} Suppose $\rho_k(W)\leq(k+1)(k-2)$. Then by Claim~\ref{Oresmall}, $\rho_k(W)=(k+1)(k-2)$.
If $G[W]$ is disconnected, say $W=W_1\cup W_2$ with $W_1\cap W_2=\emptyset$ and no edges between $W_1$
and $W_2$, then
 $\rho_{k,G}(W_1) + \rho_{k,G}(W_2) = \rho_{k,G}(W)$.
This implies $\min\{ \rho_{k,G}(W_1), \rho_{k,G}(W_2)\} \leq (k+1)(k-2)/2$, which contradicts Claim~\ref{Oresmall}.  So 
 $G[W]$ is connected, i.e., the first part of the lemma holds.

To prove the second part, choose a counter-example $G$ with the fewest vertices and let $ W\subseteq V(G)$ be a smallest   subset of $V(G)$ such that
\begin{equation}\label{j141}
 \mbox{
  $|W|\geq 2$,
 $\rho_{k,G}(W)=(k+1)(k-2)$, and $W$ does not contain a standard set.}
\end{equation}

By Fact~\ref{f1}.7, the graph $K_k$ simply does not have sets $W$ satisfying~\eqref{j141}.
So $G\neq K_k$  and thus
by Fact \ref{f2} has a separating set $\{x,y\}$.
Let $A$, $B$,  $ \widetilde G(x,y)$,
and $\check G(x,y)$ be defined as in Fact \ref{f2}.
First we show that
\begin{equation}\label{a4}
\mbox{$G[W]$ is $2$-connected.}
\end{equation}
Indeed, suppose not. Then by the first part of the lemma, $G[W]$ has a cut vertex, say $z$. Let $W_1$ and $W_2$ be two subsets of $W$ such
that  $W_1\cap W_2=\{z\}$,  $W_1\cup W_2=W$ and there are no edges between $W_1-z$ and $W_2-z$.
Then by Fact~\ref{f1}.1 and Fact~\ref{f1}.2, 
$$\rho_{k,G}(W_1)+\rho_{k,G}(W_2)=\rho_{k,G}(W)+\rho_{k,G}(\{z\})=2(k+1)(k-2).$$
So by Claim~\ref{Oresmall},  $\rho_{k,G}(W_1)=\rho_{k,G}(W_2)=(k+1)(k-2)$. Thus by the minimality of $W$, each of $W_1$ and $W_2$
contains a standard subset, a contradiction to \eqref{j141}.
 This proves~\eqref{a4}.

Let $W_A = A \cap W$ and $W_B = B \cap W$. 
Suppose $S \subseteq W_A$, $\rho_{k,G}(S) = (k+1)(k-2)$ (which implies $S \neq \emptyset$), and $y \notin S$.
Because (a) by Fact \ref{f1.1}.8, $S$ has the same potential in $\widetilde G$ as in $G$, 
(b) by Fact \ref{f2}.ii, $\widetilde G$ is also $k$-Ore, and 
(c) $G$ is a minimal counterexample, 
\begin{equation}\label{induction for finding standard sets}
\mbox {$S$ contains a standard set $W'$ in $\widetilde G(x,y)$. }
\end{equation}

We will use~\eqref{induction for finding standard sets} in Cases 1 and 3 below.


{\bf Case 1:} $W \subseteq A$. If $\{x,y\}\subseteq W$, then $\rho_{k,\widetilde G(x,y)}(W)=\rho_{k,G}(W)-2(k-1)=k(k-3)$, which
by Claim~\ref{Oresmall} means
that $W=A$. But $A$ is a standard set, a contradiction to the choice of $W$. 
So by symmetry, we may assume  that $y \notin W$.
By (\ref{induction for finding standard sets}) with $S = W$, the set $W$ contains a standard set $W'$ in $\widetilde G(x,y)$.
Since   $y\notin W'$,  by Fact \ref{f1.1}.8,
 $W'$ has the same potential in $G$ as in $\widetilde G$.
So by the minimality of $W$, $W'=W$. Furthermore,
 if $x \in W'$, then it is in the border  of $W'$ in $\widetilde G$ because $y \notin W = W'$ and $xy \in E(\widetilde G)$.
 By this and Fact \ref{f1.1}.6, we conclude that
the border   of $W'$ in $\widetilde G$ coincides with the border  of $W'=W$ in $G$.
 So $W'=W$ is also a standard set in $G$ with the same border,  a contradiction to \eqref{j141}.

{\bf Case 2:} $W\subseteq B$. 
 Let $W_B' = W - \{x,y\} + x*y$.
If $\{x,y\}\subseteq W$, then by Fact \ref{f1.1}.8 $\rho_{k,\check G(x,y)}(W'_B)=\rho_{k,G}(W)-(k+1)(k-2)=0$, which
contradicts Claim~\ref{Oresmall}.  

Suppose now that $\{x,y\}\cap W=\emptyset$. Then $W\subset V( \check G(x,y))$. By the minimality of $G$,  $W$ contains
a standard set $W'$ in $\check G(x,y)$, 
and $W'$ does not contain $x*y$. 
As in Case 1, by Fact \ref{f1.1}.8, $W'$ has the same potential in $G$ as in $\check G$, and  by Fact \ref{f1.1}.7, $W'$ has the same
 border  in $G$ as in $\check G$.
So $W'$ is also a standard set in $G$ with the same border,  a contradiction to \eqref{j141}.
Thus by the symmetry between $x$ and $y$ we may assume $\{x,y\}\cap W=\{x\}$. 
Let $i$ denote the number of neighbors of $y$ in $W$. 

{\bf Case 2A:}  $i=0$. Then $\rho_{k,\check G(x,y)}(W'_B)=\rho_{k,G}(W)=(k+1)(k-2)$,
and by the minimality of $G$, $W'_B$ contains
a standard set $W'$ in $\check G(x,y)$. If $x*y\notin W'$, then  $W'$ is standard in $G$ exactly as in the previous paragraph.
So assume $x*y\in W'$. By the case, $y$ has no neighbors in $W'$, but it does have a neighbor in $B$. This means, $x*y$ is a border
vertex of $W'$ in $\check G(x,y)$. But then the set $W'':=W'-x*y+x$ is a standard set in $G$ with the border
$W''_*=W'_*-x*y+x$,    a contradiction to \eqref{j141}.

{\bf Case 2B:}  $i\geq 1$. 
By Fact \ref{f1.1}.5, we can calculate the potential of $W'_B$ in $\check G(x,y)$ exactly:
$$\rho_{k,\check G(x,y)}(W'_B)=\rho_{k,G}(W)-2i(k-1)=k(k-3)-(i-1)2(k-1).$$ 
By Claim~\ref{Oresmall} and the definition of Ore-graphs, this yields that $W'_B=V(\check G(x,y))$ and $i=1$. 
It follows that $W=B-y$, and $y$ has exactly one neighbor, say $z$ in $W$. 
But now we are exactly in the situation described by Fact \ref{vertex to edge and back} (as $B$ is a $k$-quasi-$xy$-edge by Fact \ref{f1.1}.2), 
and so $W$ is a $k$-quasi-$xz$-vertex, and therefore a standard set,  a contradiction to \eqref{j141}.

{\bf Case 3:} $W-A\neq\emptyset$ and $W-B\neq\emptyset$. Since $\{x,y\}$ separates
$A-\{x,y\}$ from $B-\{x,y\}$ in $G$, the set $\{x,y\} \cap W$ is a separating set in $G[W]$.
Thus  by~\eqref{a4}, $\{x,y\}\subseteq W$.

{\bf Case 3A:} $W_A = A$ (equivalently, $A \subseteq W$).
We are done, since $A$ is standard by Fact \ref{f1.1}.1, Fact \ref{f1.1}.3, and the definition of a standard set.

{\bf Case 3B:} $W_B = B$ (equivalently, $B \subseteq W$).
Since $\rho_{k,G}(B)=2(k^2 - 2k - 1)$ by Fact~\ref{f1.1}.4 and $x,y$ have no common neighbors in B by \ref{f1.1}.5, we have
\begin{equation}\label{removing the second child}
\rho_{k,\widetilde G(x,y)}(W_A)=\rho_{k,G}(W)-\rho_{k,G}(B)+\rho_{k,\widetilde G(x,y)}(\{x,y\})
\end{equation}
$$=(k+1)(k-2)-2(k^2 - 2k - 1)+2(k+1)(k-2)-2(k-1)=(k+1)(k-2).$$
So for $S = W_A$ the three conditions (a), (b), (c) that imply (\ref{induction for finding standard sets})  hold,
since~\eqref{removing the second child} implies (a), and
 statements (b) and (c) are about $G$ (and hence still hold).
Thus by~(\ref{induction for finding standard sets}),
$W_A$ contains a standard set $W'$ in $\widetilde G(x,y)$.
If $|W' \cap \{x,y\}| \leq 1$, then  $W'$ is a standard set in $G$ with the same border   using the same argument as in Case 1.
If $\{x,y\}\subset W'$, then replacing in (\ref{removing the second child})   $W_A$ and $W$ with $W'$ and $ W'\cup B$ respectively, we get
$\rho_{k,G}(W' \cup B)=\rho_{k,\widetilde G(x,y)}(W')=(k+1)(k-2)$.

We claim that $W' \cup B$ has the same border  in $G$ as $W'$ does in $\widetilde G$. By Fact \ref{f1.1}.6, the only possible new elements of the border  
of $W'$ in $W' \cap A$ are $x$ or $y$.
But the only new neighbors of $x$ or $y$ are in $B$, and $B \subset W' \cup B$ so neither $x$ nor $y$ can be a \emph{new} element of the border.
And the neighbors of each $v\in B - x - y$ are in $B$, so none of the vertices in $B - x - y$ is in the border.
This means that $W' \cup B \subseteq W$ is a standard set in $G$.

{\bf Case 3C:} $W_A \neq A$ and $W_B \neq B$.
Similarly to Case 2, let $W'_B = W_B - \{x,y\} + x*y$.
By Claim \ref{Oresmall} and since each of $\check G$ and $\widetilde G$ is $k$-Ore,  $\rho_{k,\widetilde G}(W_A)\geq (k+1)(k-2)$ and $ \rho_{k,\check G}(W_B') \geq (k+1)(k-2)$.
Recall that $\{x,y\} \subset W$ by~\eqref{a4}, so $\{x,y\} =W_A\cap W_B$, and so $\rho_{k,\widetilde G}(W_A) = \rho_{k,G}(W_A) - 2(k-1)$ and by Fact~\ref{f1.1}.8, $\rho_{k,\check G}(W_B') = \rho_{k,G}(W_B) - (k+1)(k-2)$.
Therefore 
$$ \rho_{k,G}(W_A) + \rho_{k,G}(W_B) \geq (k+1)(k-2) + 2(k-1) + (k+1)(k-2) + (k+1)(k-2)=3(k+1)(k-2)+ 2(k-1) .$$
But since $W_A\cap W_B=\{x,y\}$ and $xy\notin E(G)$, by Fact~\ref{f1}.1, 
$$\rho_{k,G}(W_A)+\rho_{k,G}(W_B)=\rho_{k,G}(W)+\rho_{k,G}(\{x,y\})=3(k+1)(k-2),$$ a contradiction.
\qed

Now we will prove two statements on colorings and structure of  subgraphs not containing standard sets of $k$-Ore graphs.

\begin{lemma}\label{extr} Let  $G$ be a $k$-Ore graph.
Let $uv$ be an edge in $G$ such that
\begin{equation}\label{a61}
 \mbox{ $\rho_{k,G-uv}(W)>(k+1)(k-2)$ for every  $W\subseteq V(G-uv)$
with $2\leq |W|\leq |V(G)|-1$.}
\end{equation}
Then for each $w\in V(G)-u-v$, there is a $(k-1)$-coloring $\phi_w$ of $G-uv$ such that
$\phi_w(w)\neq \phi_w(u)=\phi_w(v)$.
\end{lemma}

{\bf Proof.} We use induction on $|V(G)|$.  For $G=K_k$, the statement is evident.
Otherwise, let $x,y,A,B, \widetilde G(x,y)$ and $\check G(x,y)$ be as in Fact~\ref{f2}.
By Fact \ref{f1.1},  $\rho_{k,G}(A)=(k+1)(k-2)$, and thus by~\eqref{a61},
  $uv \in E(G[A])$.

{\bf Case A:} $w\in A$. By the induction assumption, there exists a
$(k-1)$-coloring $\phi'_w$ of $ \widetilde G(x,y)-uv$ such that
$\phi'_w(w)\neq \phi'_w(u)=\phi'_w(v)$. 
Since $\phi'_w(x)\neq \phi'_w(y)$ and $B$ is a quasi-$xy$-edge, this coloring extends to a $(k-1)$-coloring of the whole $G-uv$.

{\bf Case B:} $w\in B-x-y$. 
Let $\phi'$ be any $(k-1)$-coloring of $ \widetilde G(x,y)-uv$. 
By Fact \ref{f2}(iv),  $xy \notin E(G)$, so $uv \neq xy$.
By the definition of $\widetilde G$, $xy \in E(\widetilde G)$, so  $\phi'(x) \neq \phi'(y)$.  
Since $ \widetilde G(x,y)$ is $k$-critical, $\phi'(u)=\phi'(v)$.

{\bf Case B1:} $\phi'(u)=\phi'(x)$. 
Let $G_0=G[B]+xw$ if $xw \notin E(G)$ and $G_0 = G[B]$ otherwise.
Then for each $W \subseteq V(G_0)$,  
\begin{equation}\label{16a61}
 \mbox{ 
$\rho_{k,G_0}(W) = \rho_{k,G}(W)$ if $\{x,w\} \not\subseteq W$ and $\rho_{k,G_0}(W) \geq \rho_{k,G}(W) - 2(k-1)$ if $\{x,w\} \subseteq W$.}
\end{equation}
Since $\{u,v\}\not\subseteq V(G_0)$,  by (\ref{a61}), $\rho_{k,G}(W) > (k+1)(k-2)$ for each $W \subseteq V(G_0)$ with $|W| > 1$.
This  together with~\eqref{16a61} imply that 
$$\rho_{k,G_0}(W)\geq \rho_{k,G}(W) - 2(k-1) >(k+1)(k-2)-2(k-1)=k(k-3)$$ 
for every $W\subseteq V(G_0)$ with $|W| > 1$.
If $|W| = 1$, then $\rho_{k,G_0}(W) = (k+1)(k-2) > k(k-3)$, and so $\rho_{k,G_0}(W)>k(k-3)$ for all nonempty $W\subseteq V(G_0)$.
By the second part of Corollary \ref{k(k-3)}, this implies that $G_0$ has a $(k-1)$-coloring $\phi''$. 
Since $G[B] \subseteq G_0[B]$, 
coloring  $\phi''$ is also a coloring of quasi-$xy$-edge $G[B]$, 
 which yields $\phi''(x)\neq \phi''(y)$.
By Fact \ref{f2}(i) and because $\phi''(x)\neq \phi''(y)$ and $\phi'(x) \neq \phi'(y)$, we can rename the colors in $\phi''$ so that $\phi''(x)=\phi'(x)$ and $\phi''(y)=\phi'(y)$,
 and obtain a $(k-1)$-coloring $\phi=\phi'\vert_{A}\cup \phi''\vert_{B}$.
By construction, $\phi(u)=\phi(x)\neq \phi(w)$.

{\bf Case B2:} $\phi'(u)\notin\{\phi'(x),\phi'(y)\}$ and $k\geq 5$. Since $B$ is a quasi-edge, $G[B]$ has a
 $(k-1)$-coloring $\phi''$ of $G[B]$ such that $\phi''(x)=\phi'(x)$ and $\phi''(y)=\phi'(y)$.
If $\phi''(w) \in \{\phi''(x),\phi''(y)\}$, then by the assumption of the case, $\phi=\phi'\vert_{A}\cup \phi''\vert_{B}$ is 
 a $(k-1)$-coloring we are looking for.
Otherwise, since $k-1\geq 4$, we can rename the colors of $\phi''$ distinct from the colors of $x$ and $y$  so that $\phi''(w)\neq \phi'(u)$
and again take $\phi=\phi'\vert_{A}\cup \phi''\vert_{B}$.

{\bf Case B3:} $\phi'(u)\notin\{\phi'(x),\phi'(y)\}$ and $k=4$. 
Let $G_0$ be obtained from $G[B]$ by adding a new vertex $z$ adjacent to $x,y$ and $w$. 
Suppose first that $G_0$ has a $3$-coloring $\phi''$.
Since $G[B] \subset G_0$ and $B$ is a quasi-$xy$-edge,  $\phi''(x) \neq \phi''(y)$.
So  the color of $z$ is distinct from $\phi''(x)$ and $\phi''(y)$, and thus because there are only $3$ colors, $\phi''(w)\in \{\phi''(x),\phi''(y)\}$. 
In this case by renaming the  colors in $\phi''$ so that $\phi''(x)=\phi'(x)$ and $\phi''(y)=\phi'(y)$, we get a required coloring of $G$. 
Now suppose that  $G_0$ has no $3$-coloring.  Then $G_0$ contains a $4$-critical subgraph $G_1$. 
Since $G_1$ is not a subgraph of $G$, it follows that $z\in V(G_1)$. 
Since $G_1$ is  $4$-critical, $\delta(G_1)\geq 4-1 = 3$, and so $\{x,y,w\}\subset V(G_1)$.
Let $W=V(G_1)$. 
Since $\rho_{4,G_0}(W)\leq 4$ by Corollary \ref{k(k-3)}, we have  $\rho_{4,G}(W-z)=\rho_{4,G}(W)-10+3(6)\leq  12$. 
So  Fact \ref{f1}.1 (because $G[A \cap W] = G[\{x,y\}] \cong 2K_1$) implies that , 
$$\rho_{4,G}(A\cup W-z)\leq \rho_{4,G}(A)+\rho_{4,G}(W-z)-2\rho_4(K_1)\leq 10+12-20=2. $$
By Claim~\ref{Oresmall}, this yields that  $A\cup W-z$ either is empty or is  $V(G)$.
But  $|A \cup W-z| \geq 3$ because $\{x,y,w\}\subset W-z$.
Also $A \cup W-z \neq V(G)$, since $G$ is $4$-Ore, and the vertex set of each $4$-Ore graph has potential $k(k-3) = 4$, since it is
$4$-extremal.   \qed

\begin{claim} \label{special vertex}
Let $G$ be a $k$-Ore graph.
Let $u$ be a vertex in $G$ such that 
\begin{equation}\label{a62}
 \mbox{$\rho_{k,G}(W)>(k+1)(k-2)$ for every  $W\subseteq V(G)-u$ with $|W|\geq 2$. }
\end{equation}
Then there exists a $(k-1)$-clique $S \subseteq V(G) - u$ such that $d_G(v) = k-1$ for all $v \in S$ and $N(S) - S$ is an independent set.
\end{claim}
{\bf Proof.}
We use induction on $|V(G)|$.
For $G=K_k$, the statement is evident.
Otherwise, let $x,y,A,B, \widetilde G(x,y)$ and $\check G(x,y)$ be as in Fact~\ref{f2}.
Then $\rho_{k,G}(A)=(k+1)(k-2)$, and so $u \in A$.

If there exists a $W \subseteq V(\check G(x,y))$ such that $|W| \geq 2$ and $\rho_{k,\check G(x,y)}(W) \leq (k+1)(k-2)$, then 
by~\eqref{a62},
$x*y \in W$. 
 So by induction, 
$\check G$   has a set $S \subseteq V(\check G(x,y))-x*y$ such that $\check G(x,y)[S] \cong K_{k-1}$,
$d_{\check G(x,y)}(v) = k-1$ for all $v \in S$,
and $N(S) - S$ is an independent set in $\check G(x,y)$.
Recall that $\check G-x*y$ is a subgraph of $G$.  So  since $u \in A$, we have 
 $S \subseteq V(G)-u$, $G[S] \cong K_{k-1}$, and $N_G(S) - S$ is an independent set in $G$.
By Fact \ref{f1.1}.7, $d_G(v) = k-1$ for all $v \in S$.
\qed

\section{Basic properties of minimal counter-examples}


The {\em closed neighborhood} of a vertex $u$ in a graph $H$ is $N_H[u] = N_H(u) \cup \{u\}$. 
We will use the following partial order on the set of graphs.
A graph $H$ is {\em smaller than} a graph $G$, if either\\
(S1) $|V(G)|>|V(H)|$, or\\
(S2)  $|V(G)|=|V(H)|$ and  $|E(G)|>|E(H)|$, or\\
(S3)  $|V(G)|=|V(H)|$,  $|E(G)|=|E(H)|$ and $G$ has fewer pairs of adjacent vertices with the same closed
neighborhood.

Note that if $H$ is a subgraph of $G$, then $H$ is smaller than $G$.
Let  $k\geq 4$ and 
$J$ be  a minimal with respect to  relation ``smaller'' counter-example to Theorem~\ref{pot theorem}:
$J$ is a $k$-critical graph with  $\rho_k(V(J))>y_k$ that is not  $k$-Ore. Since $y_k$ and the values of $\rho_k$
are always even, the restriction $\rho_k(V(J))>y_k$ is equivalent to
\begin{equation}\label{16-1}
 \rho_k(V(J))\geq y_k+2.
\end{equation}
Let $n:=|V(J)|$. 
In this section, we derive basic properties of $J$ and its colorings.

\begin{claim} \label{3-connected}
$J$ is $3$-connected.
\end{claim}
{\bf Proof.} 
Suppose that $J$ has a separating set $\{x,y\}$ and sets $A\subset V(J)$ and $B\subset V(J)$ such that $A\cap B=\{x,y\}$, 
$A\cup B=V(J)$, and no edge of $J$ connects $A-x-y$ with $B-x-y$.
By Fact~\ref{fa7} and the symmetry between $A$ and $B$, we may assume that $A$ is a $k$-quasi-$xy$-vertex and $B$ is a $k$-quasi-$xy$-edge.
It follows that the graph $\widetilde J$ obtained from $J[A]$ by inserting edge $xy$ and the graph $\check J$
 obtained from $J[B]$ by gluing $x$ with $y$ are $k$-critical.
 Then
 \begin{equation}\label{16a62}
\rho_k(V(J))\leq (\rho_k(V( \widetilde  J))+2(k-1))+(\rho_k(V(\check J))+(k+1)(k-2))-2\cdot (k+1)(k-2)
\end{equation}
$$= \rho_k(V( \widetilde J))+\rho_k(V(\check J))-k(k-3).$$
By assumption, $y_k < \rho_k(V(J))$.
By  Corollary \ref{k(k-3)}, $\rho_{k,\widetilde J}(V(\widetilde J)\leq k(k-3)$ and $\rho_k(V(\check J))\leq k(k-3)$. 
Moreover,  if $ \widetilde  J$ (respectively, $\check J$) is not a $k$-Ore graph, then by the minimality of $J$, the potential of its vertex set is at most $y_k$.
If at least one of $ \widetilde  J$ or $\check J$ is not $k$-Ore, then we get a contradiction with~\eqref{16a62}.
If both are $k$-Ore, then $J$ is $k$-Ore, which contradicts the definition of $J$.
\qed

\begin{fact}\label{aug3}
By the definition of $\rho_k$ and the assumption $\rho_k(V(J)) > y_k$, for each $v\in V(J)$,
  $$\rho_{k}(V(J)-v)=\rho_{k}(V(J))-(k+1)(k-2)+2(k-1)d(v)>$$
\begin{itemize}
	\item $y_k+k^2-3k+4$, if $d(v)=k-1  $,
	\item $y_k+k^2-k+2 $, if $ d(v)= k	 $,
	\item $y_k+k^2+k  $, if $ d(v)\geq k+1 $.
\end{itemize}
Because $y_k \geq k^2 - 5k + 2$, we see that $\rho_{k}(V(J)-v)$ is also more than 
\begin{itemize}
	\item $2k^2-8k+6 = 2(k-3)(k-1)$, if $d(v)=k-1  $,
	\item $2k^2-6k+4 = 2(k-2)(k-1)$, if $d(v)= k $,
	\item $2k^2-4k+2 = 2(k-1)^2$, if $ d(v)\geq k+1 $.
\end{itemize}
\end{fact}

Now we define graph $Y(J,R,\phi,X)$. The idea of  $Y(J,R,\phi,X)$ is that it is often smaller than $J$, and
every $(k-1)$-coloring of it extends to a $(k-1)$-coloring of $J$.

\begin{defn}\label{d1} 
For a graph $G$, a set $R\subset V(G)$ and a $(k-1)$-coloring $\phi: R\to [k-1]$ of $G[R]$, the graph $Y(G,R,\phi,X)$ is constructed as follows. 
Let $R_*$ be the border  of $R$, i.e., $R_* = \{v \in R : N(v) - R \neq \emptyset\}$.
Let $t$ be the number of colors used by $\phi$ on $R_*$. 
We may renumber the colors so that the colors used  by $\phi$ on $R_*$ are $1,\ldots,t$.
First, for $i=1,\ldots,t$, let
$R'_i$ denote the set of vertices in $V(G)-R$ adjacent in  $G$
to at least one vertex  $v\in R$ with $\phi(v)=i$.
Now, let $Y(G,R,\phi,X)$ be obtained from $G - R$ by adding a set  $X=\{x_1,\ldots,x_{t}\}$
of new vertices
such that $N(x_i) = R'_i \cup (\{x_1,\ldots,x_{t}\}-x_i)$ for $i=1,\ldots,t$.
\end{defn}

Informally, the definition  can be rephrased as follows:
For a given $R\subset V(G)$ and a $(k-1)$-coloring $\phi$ of  $G[R]$, we glue each color class of $\phi(G[R])$ into a single vertex, 
then add all possible edges between the new vertices (corresponding to the color classes) and then delete those that have no neighbors outside of $R$.
Graph $Y(G,R,\phi,X)$ will be a helpful gadget for deriving properties of $G$, since it inherits some structure from $G$.

First we will prove some  properties of $Y(J,R,\phi,X)$.

\begin{claim} \label{coloring Y}
Suppose $R \subset V(J)$ and $\phi$ is a $(k-1)$-coloring of $J[R]$.
Then $\chi(Y(J, R, \phi,X)) \geq k$.
\end{claim}
{\bf Proof.}
Let $Y = Y(J, R, \phi,X)$.
Suppose $Y$ has a $(k-1)$-coloring $\phi'$.  
By the construction of $Y$, the colors of all $x_i$ in $\phi'$ are distinct.
We can change the names of the colors so that $\phi'(x_i) = i$ for $1 \leq i \leq t$, where $t$ is given in Definition~\ref{d1}.
Again by the construction of $Y$, $\phi'(u) \neq i$ for each vertex $u \in R'_i$.
Therefore $\phi|_R \cup \phi'|_{V(J)-R}$ is a proper coloring of $J$,  a contradiction.
\qed

The next claim is a submodularity-type equation that is a direct extension of Fact \ref{f1}.1.

\begin{claim} \label{bound}
Let $R \subset V(J)$, $\phi$ be a $(k-1)$-coloring  of $J[R]$ and $Y=Y(J,R,\phi,X)$.
Let $W\subseteq V(Y)$.
If  $W\cap X=\{x_{i_1},\ldots,x_{i_q}\}$, then let
$R|_W$ denote the set of vertices $v\in R_*$ such that
$\phi(v)\in \{{i_1},\ldots,{i_q}\}$.
Then
\begin{equation}\label{j291}
\rho_{k,J}(W-X+R) = \rho_{k,Y}(W)-\rho_{k,Y}(W\cap X)+\rho_{k,J}(R)-2(k-1)|E_J(W-X,R-R|_W)|.
\end{equation}
\end{claim}
{\bf Proof.} Since $\rho_{k,J}(U)$ is a linear combination of the numbers of vertices and edges
in $J[U]$, it is enough to check that the weight of every vertex and edge of $J[W-X+R]$ is accounted
exactly once in the RHS of (\ref{j291}) and the weight of every other vertex or edge
either does not appear at all or appears once with plus and once with minus.
In particular, the weight of every  vertex and edge of $Y[W\cap X]$
appears once with plus and once with minus.\qed

By Corollary \ref{k(k-3)} and Claim~\ref{coloring Y}, $Y(J, R, \phi,X)$  contains a vertex set with potential at most $k(k-3)$.
In some instances this will not be  enough for our purposes, and we will want $Y(J, R, \phi,X)$ to contain a vertex set with potential at most $y_k$.
The next claim helps us with this.

\begin{claim} \label{potential of Y}
Let $R \subset V(J)$, $\phi$ be a $(k-1)$-coloring  of $J[R]$ and $Y=Y(J,R,\phi,X)$.
Then $Y$ contains a $k$-critical subgraph $Y'$, and so $\rho_{k, Y}(V(Y')) \leq k(k-3)$.
Furthermore, if $|R| \geq k$, then $Y$ is smaller than $J$ and \\
(a) either $Y'$ is an induced $k$-Ore subgraph, or \\
(b)  
 $\rho_{k,Y}(V(Y')) \leq y_k$.\\
Moreover, $V(Y') \cap X \neq \emptyset$.
\end{claim}

{\bf Proof.} 
By Claim \ref{coloring Y}, $Y$ has a $k$-critical subgraph $Y'$.
The bound on the potential of $V(Y')$ follows from Corollary \ref{k(k-3)}. 
In order to prove the ``Furthermore''  part, observe that if $|R| \geq k$, then $Y$ is 
smaller than $J$ by Rule (S1) in the definition of ``smaller'', since $\phi$ uses
at most $k-1<|R|$ colors on $R$. 
Because $Y'$ is a subgraph of $Y$,  $Y'$ is smaller or equal to $Y$, and so by transitivity  is smaller than $J$.
Thus, by the minimality of $J$, either $Y'$ is $k$-Ore or  $V(Y')$ 
has potential at most $y_k$.   
If $Y'$ is an induced subgraph and $k$-Ore, then (a) holds.
If $Y'$ is not induced, then by Fact \ref{f1}.10, (b) holds.
If $Y'$ is not $k$-Ore, then (b) holds by the minimality of $J$.

Since $\chi(Y')>\chi(J[R])$, $Y'$ is not a subgraph of $J$. So, $V(Y') \cap X \neq \emptyset$.
\qed

Now we will use $Y(J,R,\phi,X)$ to prove lower bounds on  potentials of nontrivial sets.

\begin{claim} \label{very small}
If $\emptyset\neq R \subsetneq V(J)$,  then $\quad\rho_{k,J}(R) \geq \rho_k(V(J)) + 2(k-1) > y_k+2(k-1)$.
\end{claim}
{\bf Proof.} 
Let $R$  be a nonempty proper subset of $V(J)$ with the smallest potential.
Since $J$ is $k$-critical, $J[R]$ has a proper coloring $\phi:R \rightarrow [k-1]$.
Let $Y=Y(J,R,\phi,X)$. 
By Claim~\ref{potential of Y}, $Y$ contains a subset $S$ with potential at most $k(k-3)$ and $S \cap X \neq \emptyset$.
Let $Z = S - X + R$.
Because $|X| \leq k-1$, by Fact \ref{f1}.7, each non-empty subset of $X$ has potential at least $(k+1)(k-2)$. So by (\ref{j291}),
\begin{equation}\label{j19}
\rho_{k,J}(Z) \leq \rho_{k,Y}(S)-\rho_{k,Y}(S\cap X)+\rho_{k,J}(R)\end{equation}
 $$\leq k(k-3)-(k+1)(k-2)+\rho_{k,J}(R)=\rho_{k,J}(R)-2(k-1).$$
Since $Z\supset R$, it is nonempty. 
So, by the minimality of the potential of $R$, we have  $Z=V(J)$.

The final statement comes from our assumption that $\rho_k(V(J)) > y_k$.
\qed

The following  fact 
 implies that $J$ contains no quasi-vertex.

\begin{claim} \label{very small potential}
For each $R \subsetneq V(J)$ with $|R| \geq 2$ and any distinct  $x,y \in R$, 
the graph $J[R]+xy$ is $(k-1)$-colorable.
\end{claim}
{\bf Proof.}
Let $R$ be a smallest subset of vertices of $J$ such that $2\leq |R|<n$ and
 for some distinct $x,y\in R$, the graph $H=J[R]+xy$ is  not $(k-1)$-colorable.
Since $J$ is $k$-critical, $xy\notin E(J)$. 

Let $H'$ be a $k$-critical subgraph  of $H$.
By the minimality of $R$, $V(H')=R$.  
By Claim~\ref{very small},  
$\rho_{k,H'}(R)=-(2k-2)+\rho_{k,J}(R) >y_k.$ 
Because $|R| < n$, 
by Rule (S1), $H'$ is smaller than $J$. So by the minimality of $J$ and because $\rho_{k,H'}(R)>y_k$, 
$H'$ is   $k$-Ore and $\rho_{k,H'}(R)=k(k-3)$. If there is an edge $e\in E(H)-E(H')$,
then $$\rho_{k,H}(R)\leq \rho_{k,H'}(R)-2(k-1)=k(k-3)-2(k-1)\leq y_k,$$
contradicting Claim~\ref{very small}.
Hence $H=H'$. Thus $H$  is $k$-Ore and
\begin{equation}\label{a72}
\rho_{k,J}(R) = \rho_{k,H}(R) + 2(k-1) = k(k-3) + 2(k-1) = (k+1)(k-2). 
\end{equation}

By Claim \ref{3-connected}, $|R_*| \geq 3$. We want to prove that
\begin{equation}\label{au1}
 \mbox{ $J[R]$ has a $(k-1)$-coloring $\psi$ such that $R_*$ is not monochromatic.}
\end{equation}

{\bf Case 1:} $\{x,y\}\subset R_*$. 
Since $|R_*| \geq 3$, we may choose $w\in R_*-x-y$.
If there exists a subset $R' \subsetneq R$ with $|R'|\geq 2$ such that $\{x,y\} \not\subset R'$ and $\rho_k(R') = (k+1)(k-2)$, then by Lemma \ref{ore4},
$H$ contains a standard set $A \subseteq R'$.
But then there exists a pair of vertices $\{a,b\} \subset A \subseteq R' \subsetneq R$ such that $J[A] + ab$
 is not $(k-1)$-colorable, which contradicts the minimality of $R$.
Otherwise, by Lemma~\ref{extr}, there is a $(k-1)$-coloring $\phi_w$ of $H-xy$ such that $\phi_w(w)\neq \phi_w(x)=\phi_w(y)$. 
Then for $\psi=\phi_w$, (\ref{au1}) holds.

{\bf Case 2:} $\{x,y\}\not\subset R_*$. 
Let $u,v$ be any vertices in $R_*$.
If $uv \in E(J)$, then (\ref{au1}) is immediately true.
Otherwise, let $H_0=J[R]+uv$. 
If $H_0$ has a $(k-1)$-coloring, then (\ref{au1}) holds. 
If not, then by the minimality of $R$, exactly as above, $H_0$ is a $k$-Ore graph. 
So, we have Case 1. 
This proves (\ref{au1}).

\medskip

Let $\psi$ satisfy (\ref{au1}).
Let $Y=Y(J,R,\psi,X)$. 
By Claim~\ref{potential of Y},  $Y$ contains a vertex set $W$ such that $\rho_{k,Y}(W) \leq k(k-3)$ and $W \cap X \neq \emptyset$.
Recall that $X$ induces a complete graph  and that by Fact~\ref{f1}.7 the 
subgraph of $K_{k-1}$ with the smallest potential is $K_1$ with  $\rho_k(V(K_1)) = (k+1)(k-2)$ and the subgraph with the second smallest potential is $K_{k-1}$ with  $\rho_k(V(K_{k-1})) = 2(k-2)(k-1)$.
 This together with~\eqref{a72} and the choice of $W$ yields
\begin{equation}\label{au2}
  \rho_{k,J}(W - X + R) \leq \rho_{k,Y}(W)-\rho_{k,Y}(X \cap W)+\rho_{k,J}(R)\leq  k(k-3)-(k+1)(k-2)+(k+1)(k-2)=k(k-3).
\end{equation}
Since $W-X+R\supset R$, we have $|W-X+R|\geq 2$. 
By Fact \ref{f1}.8, $y_k+(2k-2) \geq k(k-3)$. This together  with Claim~\ref{very small} yields $W-X+R=V(J)$. 
If $|W\cap X|\geq 2$, then we get the stronger bound $\rho_k(X\cap W)\geq 2(k-1)(k-2)$, and so  the inequality in (\ref{au2}) strengthens to 
 $$ \rho_{k,J}(W - X + R) \leq k(k-3)-2(k-1)(k-2) +(k-2)(k+1)=2k-6\leq y_k,$$ a contradiction.
Thus $|X \cap W| = 1$.
Because $R_*$ is not monochromatic and $|X \cap W| = 1$, there is a vertex $z\in R_*-W$. Then by
(\ref{j291}), instead of (\ref{au2}) we have
$$ \rho_{k,J}(W - X + R) \leq  k(k-3)-(k+1)(k-2) +(k+1)(k-2)-2k+2=k^2-5k+2\leq y_k,$$
 a contradiction.
\qed

\begin{claim} \label{one cluster}
Let $X$ be a $(k-1)$-clique in $J$, $u,v \in X$, $N(u) - X = \{a\}$, and $N(v) - X = \{b\}$.
Then $a = b$.
\end{claim}
{\bf Proof.} 
Assume $a \neq b$.
Let $J'= J - u - v + ab$ if $ab \notin E(J)$ and $J' = J - u - v$ otherwise.
By Claim~\ref{very small potential},  $J'$ has a $(k-1)$-coloring $\phi$. 
Because $d(u) = d(v) = k-1$, each of the sets $C_a = \{1,\ldots,k-1\} - \cup_{w \in N(u)- v}\phi(w)$ and 
$C_b = \{1,\ldots,k-1\} - \cup_{w \in N(v)- u}\phi(w)$  is nonempty.
Since $\phi(a)\neq\phi(b)$ and $(N(u) - a) = (N(v) - b)$, the sets $C_a$ and $C_b$ are different.
Therefore $\phi$ can be extended to $u$ and $v$.
But then we get a $(k-1)$-coloring of $J$,  a contradiction.
\qed

\begin{claim} \label{minus edge} 
$J$ does not contain $K_k-e$.
\end{claim}
{\bf Proof.} 
Suppose $J[R]=K_k-e$. 
The only $k$-critical graph on $k$ vertices is $K_k$, which is $k$-Ore.
By assumption $J$ is not $k$-Ore, so $R\neq V(J)$, but adding the missing
edge to $J[R]$ creates a $k$-chromatic graph on $R$, a contradiction to
Claim~\ref{very small potential}.
\qed

\section{Clusters and sets with small potential}
Recall that in Section $3$ we defined a relation ``smaller,'' and 
 $J$  is a ``smallest'' counterexample to Theorem \ref{pot theorem}: it is a $k$-critical 
	graph with
	 $\rho_k(V(J)) > y_k$ and is not $k$-Ore.

\begin{defn} For $S\subseteq V(J)$, an
 $S$-{\em cluster} is an inclusion maximal set $R \subseteq S$ such that for every $x \in R$, $d(x) = k-1$ and for every  $x,y \in R$, $N[x] = N[y]$.
A {\em cluster } is a $V(J)$-cluster.
\end{defn}

Each cluster is a clique of vertices with degree $k-1$.
In this section, results on clusters will help us to derive the main lower bound on the potentials of nontrivial vertex sets,   Lemma \ref{small potential},
which in turn will help us to prove stronger results on the structure of clusters in $J$.

Having the same closed neighborhood is an equivalence relation, and so the set of clusters is a partition of the 
set  of the vertices with degree $k-1$. Thus the following fact holds.

\begin{fact}\label{clusters partition}
Every vertex with degree $k-1$ is in a unique cluster. 
\end{fact}

Furthermore, if the only $S$-cluster  is the empty set, then every vertex in $S$ has degree at least $k$.
By definition, if a cluster $T$ is contained in a vertex set $S$, then $T$ is also an $S$-cluster.



\begin{claim} \label{unique cluster} Every cluster $T$ satisfies $|T| \leq k-3$. Furthermore, for every
 $(k-1)$-clique $X$ in $J$, (i)~there is a unique $X$-cluster $T$  (possibly $T = \emptyset$), and (ii)~every non-empty $X$-cluster 
 is a cluster (in other words, every cluster is either contained in $X$ or disjoint from $X$).
Each $(k-1)$-clique in $J$  contains at least $2$ vertices of degree at least $k$.
\end{claim}
{\bf Proof.} If $T$ is a cluster with $|T| \geq k-2$, then $T \cup N(T) \supseteq K_k - e$, a contradiction to Claim~\ref{minus edge}.

Let $X$ be a $(k-1)$-clique  in $J$.
Two distinct $X$-clusters  would contradict Claim \ref{one cluster}. If $T$ is a non-empty $X$-cluster contained in a larger cluster $T'$,
then each $v\in T'-X$ has to be adjacent to each vertex in $X$, and so $J$ contains clique $X\cup T'$ of size at least $k$, a contradiction. 

The final statement is proven as follows: by Fact \ref{clusters partition} only vertices  in  clusters  have degree $k-1$, 
by (i) $X$ contains at most one cluster which, if exists, has  at most $k-3$ of the $k-1$ vertices in $X$ by the first part of this claim.
\qed

\begin{claim} \label{standard sets in J}
Let $R \subset V(J)$, let $\phi$ be a $(k-1)$-coloring of $J[R]$, and let $Y = Y(J,R,\phi,X)$ be as in Definition \ref{d1}.
If $Y$ is $k$-Ore, then $Y \cong K_{k}$.
\end{claim}
{\bf Proof.}
Suppose that $Y$ is a $k$-Ore graph distinct from $K_k$.
Let a separating set $\{x,y\}$, vertex subsets $A=A(Y,x,y)$ and $B=B(Y,x,y)$, and graphs $\widetilde Y(x,y)$ and $\check Y(x,y)$ be as in Fact~\ref{f2}.
Since $Y[X]$ is a clique and $E_{Y}(A-x-y, B-x-y) = \emptyset$, either $X \subseteq A$ or $X \subseteq B $.
Since $xy\notin E(Y)$ we may assume that either $X\subset A-y$ or $X\subset B-y$.
By construction, $Y-J=X$, so this implies that $B-x$ or $A-x$, respectively, is a subgraph of $J$.
We will show that this
 is impossible.

Suppose first that $X\subseteq A-y$. The graph $\check Y-x*y$ is a subgraph of $J$, namely, it is $J[B-x-y]$,
and by Fact \ref{f1.1}.7
\begin{equation}\label{aug4'}
\mbox{$d_{\check Y}(v)=d_J(v)$ for every $v\in B-x-y$.}
\end{equation}
If $\check Y-x*y$ has a vertex subset $S$ with $|S|\geq 2$ of potential at most $(k+1)(k-2)$, then by Lemma~\ref{ore4}, $S$ contains a standard set $S'$.
But  each standard set $S'$ has two vertices $u$ and $w$ such that $Y[S']+uw$ is not $(k-1)$-colorable.
This contradicts Claim \ref{very small potential}.
Thus $\rho_{k,\check Y}(S)>(k+1)(k-2)$ for every $S\subseteq V(\check Y)-x*y$ with $|S|\geq 2$.
Then by Claim~\ref{special vertex}, there exists an $S \subseteq V(\check Y)-x*y=B-x-y$ such that $\check Y[S] \cong K_{k-1}$,
 and $d_{\check Y}(v) = k-1$ for all $v \in S$.
By (\ref{aug4'}), this contradicts the last part of Claim \ref{unique cluster}.

Now suppose that $X \subseteq B-y$. 
The graph $ \widetilde Y -x$ is a subgraph of $J$, namely, it is $J[A-x]$,
and similarly to (\ref{aug4'}), by Fact \ref{f1.1}.6
\begin{equation}\label{aug5}
\mbox{$d_{ \widetilde Y}(v)=d_J(v)$ for every $v\in A-x-y$.}
\end{equation}
As above, $\rho_{k, \widetilde Y }(S)>(k+1)(k-2)$ for every $S\subseteq V( \widetilde Y )-x$ with $|S|\geq 2$.
So again by Claim \ref{special vertex}, there exists an $S' \subseteq V( \widetilde Y )-x=A-x$ such that $ \widetilde Y [S'] \cong K_{k-1}$, 
and $d_{ \widetilde Y }(v) = k-1$ for all $v \in S'$.
But $|S'- y| \geq k-2$, which together with (\ref{aug5}) contradicts Claim \ref{unique cluster}.
\qed

\begin{claim} \label{k edges}
For every partition $(A,B)$ of $V(J)$ with $2\leq |A|\leq n-2$,
$|E_J(A,B)|\geq k$.
\end{claim}
{\bf Proof.} Let $A_*$ (respectively, $B_*$) be the set of vertices in $A$ (respectively, $B$)
 that have neighbors in $B$ (respectively, $A$). Since $J$ is $3$-connected, $|A_*|\geq 3$ and
 $|B_*|\geq 3$. So by Claim~\ref{very small potential}, $J[A]$ has a $(k-1)$-coloring $\phi_A$
 such that $A_*$ is not monochromatic, and  $J[B]$ has a $(k-1)$-coloring $\phi_B$
 such that $B_*$ is not monochromatic. But Gallai and Toft~(see~\cite[p. 157]{Toft2})
 independently proved that if $|E_J(A,B)|\leq k-1$, then either $A_*$ is monochromatic in
 every $(k-1)$-coloring of $J[A]$ or $B_*$ is monochromatic in
 every $(k-1)$-coloring of $J[B]$.
 So,  $|E_J(A,B)|\geq k$.\qed

Sometimes below, our goal will be to extend to $J$ a coloring $\phi$ of $J[R]$ for some $R$ and $\phi$.
Recall that $Y(J,R,\phi,X)$ is obtained from $J$ by replacing the vertices of $R$ with a clique whose vertices are 
the color classes of $\phi$ with at least one element in the border  of $R$ (which we called $R_*$).
One of the ways we will control $\phi$ is to add edge(s) to $R$ before we generate a $(k-1)$-coloring $\phi$ using Claim \ref{very small potential} 
and a lemma below.
Our next lemma describes how edges can be placed in $R$ so that no  color class of $\phi$ is too large.
The proof of this lemma will use the following old result of Hakimi.

\begin{theorem}[Hakimi~\cite{Hak}]\label{thak} Let $(w_1,\ldots,w_s)$ be a list of nonnegative integers with $w_1\geq \ldots\geq w_s$.
Then there is a loopless multigraph $F$ with vertex set $\{u_1,\ldots,u_s\}$ such that $d_F(u_j)=w_j$ for all $j=1,\ldots,s$ if and only if
$z=w_1+\ldots+w_s$ is even and $w_1\leq w_2+\ldots+w_s$.
\end{theorem}

For technical reasons,
in one specific case of the lemma below we will allow for a hyperedge of size $3$. Recall that
an {\em independent set} in a hypergraph is a set that contains no edge.

\begin{lemma}\label{lem1 - k} Let $i'\geq 1$ and $ s\geq 2$ be  integers.
Let $R_*=\{u_1,\ldots,u_s\}$ be a vertex set. Then for each $z\geq 2i'$ and any
integral positive weight function  $w\,:\,R_* \to \{1,2,\ldots\}$  such that $w(u_1)+\ldots+w(u_s)=z$ and $w(u_1)\geq w(u_2)\geq\ldots\geq w(u_s)$,
 there exists a graph $ H$ with $V( H)=R_*$ and $|E( H)|\leq i'$ such that for each $1\leq j\leq s$, $d_ H(u_j)\leq w(u_j)$, and for every independent set $M$ in 
 $ H$ with $|M|\geq 2$,
\begin{equation}\label{I01}
\mbox{ $\sum_{u\in R_*-M}w(u)\geq i'$.}
\end{equation}
Moreover, if  $s\geq 3$, $i'\geq 1$, and $z>2i'$, then at least one of the three stronger statements below holds:\\
(i) such $ H$ with Property (\ref{I01}) could be chosen  as a graph with at most $i'-1$ edges,  or\\
(ii) such $ H$ with Property (\ref{I01}) could be chosen as a hypergraph instead of a graph with at most $i'-1$ graph edges and one edge of size $3$, or \\
(iii)  the weight arrangement is {\em $i'$-special}, which means that $s=i'+1$
and $w(u_2)=\ldots=w(u_s)=1$.
\end{lemma}
{\bf Proof.}
It is easy to check that the main part of the statement, and (ii) of the ``moreover'' part hold for $i'=1$.
So assume $i' \geq 2$.

{\bf Case 1:} $w(u_2)+\ldots+w(u_s)\leq i'-1$. We make $E( H)=\{u_1u_j\,:\; 2\leq j\leq s\}$.
If $M$ is any independent set with $|M|\geq 2$,
then $u_1\notin M$ and $w(u_1)\geq 2i'-(i'-1)$ yielding (\ref{I01}).
To prove the ``Moreover'' part in this case, observe that  our $ H$ has at most $ i'-1$ edges.

{\bf Case 2:} $s\geq 2i'+1$. Let the edge set of $H$ consist of the matching $\{u_1u_2,\ldots,u_{2i'-1}u_{2i'}\}$. 
Every independent set misses at least one end of each edge in $H$, which implies~(\ref{I01}).
Moreover, if  $z\geq 2i'+1$, then we extend
edge $u_1u_2$ to the hyperedge $\{u_1,u_2,u_{2i'+1}\}$. Then Claim~(ii) of the ``Moreover'' part of the lemma holds.

{\bf Case 3:} $w(u_2)+\ldots+w(u_s)\geq i'$ and $s\leq 2i'$.  
 Since $s\leq 2i'$, 
there exists an auxiliary integral weight function $w'\,:\,R_* \to \{1,2,\ldots\}$  such that 
\begin{equation}\label{1010}
 \mbox{$w'(u_2)+\ldots+w'(u_s)\geq i'$,
$w'(u_1)+\ldots+w'(u_s)=2i'$ and $w'(u_j)\leq w(u_j)$  $\forall\, j=1,\ldots,s$.}
\end{equation}
 By~\eqref{1010} and Theorem~\ref{thak}, there exists a 
loopless multigraph $ H'$ with vertex set $\{u_1,\ldots,u_s\}$ such that 
$d_{H'}(u_j) = w'(u_j)$ for all $j$.
We obtain a graph $ H$ from the multigraph $ H'$ by replacing each set of multiple edges with a single edge.
Every independent set in $ H$ is also independent in $ H'$.
For every independent set $M$ in $ H'$, each of its $i'$ edges has an end outside of $M$, so
$$\sum_{u\in R_*-M} w(u)\geq \sum_{u\in R_*-M} w'(u)=\sum_{u\in R_*-M} d_{ H'}(u)\geq |E( H')|=i'.$$
 This yields (\ref{I01}).
Note that in this case, (\ref{I01}) holds for {\em every} independent set $M$,
even if $|M|=1$.

Now we prove the ``Moreover'' part of the statement.
If $ H'$ had any multiple edge, then we satisfy (i) and are done. Suppose, $ H'$ is simple.
Since $z>2i'$,  $w'(u_\ell)<w(u_\ell)$ for some $1\leq \ell\leq s$. 
If $ H-u_\ell$ has an edge $e$, then after enlarging $e$ to $e+u_l$ we still keep (\ref{I01}). 
This instance satisfies (ii), and we are done.
Otherwise $u_\ell$ is incident to every edge of $ H= H'$, and so $ H$ is a star with center $u_\ell$ and $i'\geq 2$ edges.
Each such star has only one central vertex, so
  every other vertex $u_j$ satisfies $w(u_j) =w'(u_j)= d_ H(u_j) = 1$.
By definition, this means that the weight arrangement is $i'$-special.
So we satisfy (iii) and are done.
\qed

Recall that $\rho_{k, K_{k-1}}(V(K_{k-1})) = 2(k-1)(k-2)$.
Importantly, this is larger than the potential of a standard set. Our main lower bound on the potentials
of nontrivial vertex sets is the following.

\begin{lemma} \label{small potential}
If $R \subsetneq V(J)$ and $2\leq |R|\leq n-2$, then $\rho_k(R) \geq 2(k-1)(k-2)$.
Moreover, if $\rho_k(R) = 2(k-1)(k-2)$, then  $J[R]=K_{k-1}$.
\end{lemma}
{\bf Proof.}
Assume that the lemma does not hold.
Let $i$ be the smallest integer such that there exists $R \subsetneq V(J)$
with 
\begin{equation}\label{j14'1}
 2\leq |R|\leq n-2, \, J[R]\neq K_{k-1},\,\rho_k(R) \leq 2(k-1)(k-2),
\end{equation}
 and
\begin{equation}\label{j14'}
y_k+2i(k-1)<\rho_k(R)\leq y_k+2(i+1)(k-1).
\end{equation}
It is  important  that we are only minimizing $i$, and not necessarily minimizing $\rho_k(R)$.
By Claim~\ref{very small}, $i\geq 1$. In order for both~\eqref{j14'1} and~\eqref{j14'} to hold,
since
$$y_k+(k+1)(k-1)\geq k^2-5k+2+(k+1)(k-1)>2(k-1)(k-2),$$
 $i\leq \frac{k}{2}$.
By the integrality, 
 $i\leq \left\lfloor\frac{k}{2}\right\rfloor$.
Moreover, if $k=4$ then $y_k=\max\{2\cdot 4-6,4^2-5\cdot 4+2\}=2$ and
so $y_4+4(4-1)=14>12=2(4-1)(4-2)$. Thus
\begin{equation}\label{se1}
i\leq \left\lfloor\frac{k}{2}\right\rfloor,\;\mbox{moreover,  if $k=4$ then $i=1$.}
\end{equation}

Let $R$ be a smallest subset of $ V(J)$
for which~\eqref{j14'1} and~(\ref{j14'}) hold. 
By Fact \ref{f1}.7,  each vertex set $S$ with $2\leq |S|\leq k-1$ has potential at least $2(k-1)(k-2)$,  with equality only when $S$
 induces a $K_{k-1}$.
So~\eqref{j14'1} yields  $|R|\geq k$. 
Thus by Claim \ref{potential of Y}, for any proper $(k-1)$-coloring $\phi$ of $J[R]$,  graph $Y(J,R,\phi,X)$ is smaller than $J$.

Let $Q = V(J) - R$, and for $u\in R$, let $w(u)=|N(u)\cap Q|$.
By Definition \ref{def1}, $R_* = \{u \in R\,:\, w(u)\geq 1\}$.
Let $R_*=\{u_1,\ldots,u_s\}$ and $w(u_1)\geq \ldots\geq w(u_s)$.
By Claim~\ref{k edges}, 
\begin{equation}\label{edges to add}
z:=\sum_{i=1}^sw(u_i)=|E_J(R,V(J)-R)|\geq k.
\end{equation}
By Claim~\ref{3-connected}, $s \geq 3$.

We will consider four cases, and the first is the main one. 

{\bf Case 1:} There is a $(k-1)$-coloring $\phi$ of $J[R]$ such that for every color class $C$ of $\phi$ with $|C \cap R_*| \geq 2$ either
\begin{equation}\label{low weight color classes}
\sum\nolimits_{u \in R_* - C} w(u) \geq i 
\end{equation}
or
\begin{equation}\label{low weight color classes - case 4}
\mbox{$\sum_{u \in R_* - C} w(u) = i - 1$ and $  \sum_{u \in C} w(u) \leq k-2 . $}
\end{equation}
Let $Y=Y(J,R,\phi,X)$ be as in Definition \ref{d1}.

Let $Y'$ be the $k$-critical subgraph of $Y$ described in Claim~\ref{potential of Y}. 
Let $W=V(Y')$ and
 $X' = X \cap W$.
Since $|R|\geq k$, by Claim~\ref{potential of Y}, $X' \neq \emptyset$ and one of the following 
two statements is  true: (a) $Y[W]$ contains a $k$-Ore graph, or (b) $\rho_{k,Y}(W) \leq y_k$.
Because $X' \neq \emptyset$, by Fact \ref{f1}.7, $\rho_{k,Y}(X') \geq (k+1)(k-2)$.
By Corollary \ref{k(k-3)}, $\rho_{k,Y}(W) \leq \rho_{k,Y'}(W) \leq k(k-3)$.
By (\ref{j291}) and our bounds on  $\rho_{k,Y}(X')$ and $\rho_{k,Y}(W)$,
\begin{equation}\label{aug1'}
  \rho_{k,J}(W - X + R) \leq \rho_{k,Y}(W) - \rho_{k,Y}(X') + \rho_{k,J}(R)\leq \rho_{k,J}(R)-2(k-1).
\end{equation}
Because $|W-X+R| \geq |R| \geq k$ (which implies $J[W-X+R] \neq K_{k-1}$) and the choice of $i$, $|W-X+R| \geq n-1$.
Suppose first that $|W-X+R|= n-1$.
By Fact \ref{aug3}, $\rho_{k,J}(W - X + R) >  y_k+k^2-3k+4$. By~\eqref{aug1'}, $\rho_{k,J}(R) > y_k + k^2 - k + 2 \geq 2(k-1)(k-2)$, contradicting the choice of $R$.
So 
\begin{equation}\label{june24}
W-X+R = V(J).
\end{equation}

We claim that $Y$ is a $k$-clique and
 will prove this in three steps.
We will show, in order, that
\begin{equation}\label{june242}
\mbox{ (A) $|X'| \geq 2$, (B) $Y'$ is an induced $k$-Ore graph, and (C) $Y' = Y$.}
 \end{equation}
These  three steps prove that $Y$ is $k$-Ore; Claim \ref{standard sets in J} yields that if $Y$ is $k$-Ore, then $Y$ is a $k$-clique.
 
Suppose $X'=\{x_j\}$.
Then $W = V(Y') - X+x_j$.
Let $R_j = \{u \in R_* : \phi(u) = c_j\}$, where $c_j$ is the $j$th color class in coloring $\phi$.
 If $|R_j|=1$, then $Y' \cong J[W - x_j \cup R_j]$, which is a subgraph of $J$.
Because $|R| \geq k > 1 = |R_j|$, $Y'$ is a proper subgraph of $J$, but $k$-critical graphs do not have $k$-chromatic proper subgraphs.
Thus $|R_j|\geq 2$. If~\eqref{low weight color classes - case 4} holds for $R_j$, then $d_{Y'}(x_j) \leq \sum_{u \in R_j} w(u) \leq k-2$, but
the $k$-critical graph $Y'$ cannot have vertices of degree less than $k-1$. Otherwise, by ~\eqref{low weight color classes}, 
 at least $i$ edges connect the vertices in $R_*-R_j$ with $Q$.
Those edges are a part of $E_J(W-X,R-R|_W)$ in (\ref{j291}), an edge set we conservatively estimated to be empty in the calculation of (\ref{aug1'}).
Adjusting that calculation to account for our new bound on $E_J(W-X,R-R|_W)$, instead of (\ref{aug1'}) we get 
$$	\rho_{k,J}(W - \{x_j\} + R)  \leq  k(k-3)-(k-2)(k+1) + \rho_{k,J}(R) - 2i(k-1) = \rho_{k,J}(R)-2(i+1)(k-1).$$
By (\ref{j14'}), we have  $\rho_{k,J}(R)-2(i+1)(k-1) \leq y_k$, which together with $i \geq 0$,~\eqref{aug1'}, and~\eqref{june24} imply $\rho_{k,J}(V(J)) \leq y_k$, contradicting our assumption that $\rho_{k,J}(V(J)) > y_k$.
This proves (A).

If $Y'$ is not an induced $k$-Ore graph, then by Claim \ref{potential of Y}, $\rho_{k,Y}(W) \leq y_k$.
We adjust the calculation of (\ref{aug1'}) from (\ref{j291}) again:  we use the inequality  $\rho_{k,Y}(W)\leq y_k$ 
and that part (A) implies $2\leq |X'|\leq k-1$, so $X'$ has potential at least $2(k-1)(k-2)$ by Fact \ref{f1}.7. The result is 
$$\rho_{k,J}(W-X+R) \leq y_k-2(k-1)(k-2)+ \rho_{k,J}(R).$$
We apply \eqref{june24} to the left hand side and \eqref{j14'1} to the right hand side to see that $\rho_{k,J}(V(J))\leq y_k$, which contradicts our choice of $J$.
This proves (B). 

By (\ref{june24}) and (B), 
 if $Y' \neq Y$, then $X' \neq X$. 
If $X' \neq X$, then we again get a stronger bound on $E_J(W-X,R-R|_W)$ in (\ref{j291}) (this time it contains at least one element), and so we may again adjust the bound in (\ref{aug1'}) by subtracting $2(k-1)$. 
We augment the application of (\ref{j291}) to (\ref{aug1'}) in two ways: the extra  $-2(k-1)$ described in the previous sentence and that $|X'| \geq 2$ implies $\rho_{k,Y}(X') \geq 2(k-1)(k-2)$ by Fact \ref{f1}.7; those augmentations and \eqref{j14'1} produce 
$$\rho_{k,J}(W - X + R) \leq -2(k-1)+ k(k-3) - \rho_{k,Y}(X') + \rho_{k,J}(R) \leq k^2-5k+2\leq y_k,$$
which contradicts our choice of $J$ for the same reason as the proof to part (B).
This proves (C); and therefore $Y = Y(J,R,\phi,X)$ is $k$-Ore.
By Claim \ref{standard sets in J}, $Y = K_k$.

Recall that  $\{Q,R\}$ is a  partition of $V(J)$.
Then by \eqref{june24} and (\ref{june242}C), $\{Q,X\}$ is a  partition of $V(Y)$.
Let $t=|X|$ (and $t=|X'|$ also by (\ref{june242}C)).
By~(\ref{june242}A),  $t \geq 2$.
Because $|R| \leq n-2$,  $|Q| \geq 2$. So since $V(Y) = X \cup Q$, we have $t \leq k-2$. Then $J$ is obtained from $J[R]$ by 
adding $k-t$ vertices and at least $\binom{k}{2}-\binom{t}{2}$ edges (a vertex in $Q$ may be adjacent to more than one vertex in a color
class of $\phi$). So
\begin{equation}\label{june243}
\rho_k(V(J))\leq \rho_k(R)+(k-t)(k+1)(k-2)-\left(\binom{k}{2}-\binom{t}{2}\right)2(k-1).
 \end{equation}

We know that $2\leq t\leq k-2$.
We will show that $t=2$. This  is clear when $k=4$. Let $k\geq 5$.
Denote the RHS of~\eqref{june243} by $\mu(k,t,R)$. 
For fixed $k$ and $R$, $\mu(k,t,R)$ is quadratic in $t$
with a positive  coefficient at $t^2$.
So, if $3\leq t\leq k-2$, then
 $\mu(k,t,R)\leq \max\{\mu(k,3,R),\mu(k,k-2,R)\}$. Furthermore, 
 $$\mu(k,k-2,R)=\rho_k(R)+2(k+1)(k-2)-k(k-1)^2+(k-1)(k-2)(k-3)$$
 $$\leq 2(k-1)(k-2)-2k^2+8k-10=2k-6\leq y_k,
 $$
and  when  $k\geq 5$,
 $$\mu(k,3,R)=\rho_k(R)+(k-3)(k+1)(k-2)-k(k-1)^2+6(k-1)\leq 2(k-1)(k-2)-2k^2+6k=4\leq y_k.
 $$
Since $\rho_k(V(J))>y_k$ by the choice of $J$, we conclude that $t=2$ and $J[Q]=K_{k-2}$. Moreover, $\mu(4,2,R)\leq 2(4-1)(4-2)+20-(6-1)6=2=y_4$, so $k\geq 5$.
Thus by~\eqref{june243},
\begin{equation}\label{june245}
\rho_k(V(J))\leq \mu(k,2,R)=\rho_k(R)+ (k-2)^2(k+1)-k(k-1)^2+2(k-1)=\rho_k(R)-k^2+k+2. 
 \end{equation}
Hence by~\eqref{june245} and~\eqref{16-1}, 
$$\rho_k(R)\geq \rho_k(V(J))+k^2-k-2\geq y_k+k^2-k.
$$
Plugging in the values of $y_k$, we get
\begin{equation}\label{june244}
\mbox{ 
$\rho_k(R)\geq 2(k-1)(k-2)-2$ for $k\geq 6$ and $\rho_5(R)\geq 2(5-1)(5-2)=24$.}
 \end{equation}
 Since for $k\geq 5$, $2(k-1)(k-2)-2>y_k+4(k-1)$, we have 
\begin{equation}\label{i geq 2}
i\geq 2.
\end{equation}
Also we  conclude that each $v\in Q$ has exactly two neighbors in $R$, since otherwise the upper bound on
$\rho_k(V(J))$ in~\eqref{june243} and~\eqref{june245} would be stronger by $2(k-1)$, and this combined with~\eqref{j14'1} would lead to
$$\rho_k(V(J))\leq   -2(k-1)+2(k-1)(k-2)-k^2+k+2=k^2-7k+8\leq y_k.$$



Let $Q = \{v_1, \ldots, v_{k-2}\}$ and let $N(v_j) \cap R = \{u_{j,1}, u_{j,2}\}$ for $j=1,\ldots,k-2$.
If there exists a proper $(k-1)$-coloring $\phi'$ of $J[R]$ such that $\phi'(u_{j,1}) = \phi'(u_{j,2})$ for some $j$, then $\phi'$
may be extended to all of $J$ greedily by first coloring $Q - v_j$ and at the end coloring $v_j$ (at each step at most $k-2$ colors must be avoided).
Similarly,  if   $\{\phi'(u_{j,1}), \phi'(u_{j,2})\} \neq \{\phi'(u_{j',1}), \phi'(u_{j',2})\}$ for some $j\neq j'$, then
 $\phi'$ may be extended to all of $J$ greedily by first coloring $X - v_j - v_{j'}$ and at the end coloring $v_j$ and $v_{j'}$. 
 Thus for any proper $(k-1)$-coloring $\phi'$ of $J[R]$,
 \begin{equation}\label{june246}
\mbox{for all $1\leq j,j'\leq k-2$, $\quad\phi'(u_{j,1}) \neq \phi'(u_{j,2})$ and 
$\{\phi'(u_{j,1}), \phi'(u_{j,2})\} = \{\phi'(u_{j',1}), \phi'(u_{j',2})\}$.}
 \end{equation}
  
  Because $3 \leq s = |R_*|$ by Claim~\ref{3-connected}, there exist distinct vertices $v',v'' \in Q$ such that $N(v') \cap R \neq N(v'') \cap R$. By symmetry, we may  assume $u_{1,1} \notin N(v_2)$.
Let $J^*$ be obtained from $J[R]$ by adding edges $e_1 = u_{1,1}u_{2,1}$ and $e_2 = u_{1,1}u_{2,2}$ when they do not exist.
By~\eqref{june246},  $\chi(J^*) \geq k$, which implies that $e_1$ or $e_2$ was not in $J$. 
Thus
$J^*$ contains a $k$-critical subgraph $J^\circ$.  By the minimality of $J$ and the fact that $|J^\circ|\leq |R|\leq |J|-2$, graph $J^\circ$ is $k$-Ore or $\rho_{k,J^*}(V(J^\circ)) \leq y_k$.
Because $J^\circ = J[V(J^\circ)] + S$ for some nonempty $S \subseteq \{e_1, e_2\}$, we know that $J[V(J^\circ)]$ is not a clique.
If $V(J^\circ)$ satisfies~\eqref{j14'1}, then the minimality of $i$ applied to~\eqref{j14'} states that $\rho_{k,J}(V(J^\circ)) > y_k+2i(k-1)$; 
if $V(J^\circ)$ does not satisfy~\eqref{j14'1}, then (as $\chi(J^\circ) = k$ implies $|V(J^\circ)| \geq k$, and the other conditions are proved above) $$\rho_{k,J}(V(J^\circ)) > 2(k-1)(k-2) \geq y_k+2i(k-1).$$ 
Since $i\geq 2$ by~\eqref{i geq 2} and because we have added at most two edges, by~\eqref{j14'},
 $\rho_{k,J^*}(V(J^\circ))\geq \rho_{k,J}(V(J^\circ)) -4(k-1)> y_k$, and so $J^\circ$ is $k$-Ore.
 Recall by Fact \ref{f1.1}.9, the vertex sets of $k$-Ore graphs have potential $k(k-3)$. So, $\rho_{k,J^*}(V(J^\circ)) \leq k(k-3)$ and our above inequality becomes by Fact \ref{f1}.8
$$\rho_{k,J}(V(J^\circ))\leq \rho_{k,J^*}(V(J^\circ)) + 4(k-1) \leq k(k-3)+4(k-1)\leq y_k+3(2(k-1)).$$
 Hence $V(J^\circ)$ satisfies~\eqref{j14'} for some $i\leq 2$. By~\eqref{i geq 2}, we have $i=2$.
By the minimality of $i$ and of $|R|$, this gives
 \begin{equation}\label{june25}
\mbox{$i=2$, $V(J^\circ)=R$ and $J[R]=J^\circ-e_1-e_2$.}
 \end{equation}
Also, since $i = 2$ and $(k+1)(k-2)\leq y_k+2(2(k-1))$, the minimality of the potential of $R$ in (\ref{j14'}), and Fact \ref{f1}.(2,3,7) imply
\begin{equation}\label{june252}
\mbox{$J[R]$ contains no set with potential at most $(k+1)(k-2)$.}
 \end{equation}
 
For all $S \subseteq V(J^\circ) - u_{1,1}$, we have  $J^*[S] \cong J[S]$. Thus, by~\eqref{june25} and~\eqref{june252},  Claim \ref{special vertex}
applies to $J^\circ=J^*$ and $u_{1,1}$. By this claim, $J^\circ - u_{1,1}$ contains a clique $S$ of order $k-1$ such that each vertex in $S$ has degree $k-1$.  
Since $u_{1,1}\notin S$, $S$ is also a clique in $J$.
Since $e_1, e_2 \subset N(Q) =R_*$,  if $u \in S - N(Q)$ then $d_J(u) = k-1$.
Because $N(Q)$ is $2$-colorable (by~\eqref{june246}), this implies that there is an $S' \subset S$ with $|S'| \geq k-3$ such that $d_J(u) = k-1$ for all $u \in S'$.
Each vertex of $S'$ is in a cluster by Fact \ref{clusters partition}, and Claim \ref{unique cluster} says that all of $S'$ is one cluster and that $|S'| = k-3$.
Let $\{u'\} = N(S') - S$. Then by Fact \ref{f1} (parts 1, 2, and 4)
$$\rho_{k,J}(S \cup u') \leq \rho_{k,J}(S) + \rho_{k,J}(\{u'\}) - 2(k-1)(k-3) = k^2 + k - 4<y_k+3(2(k-1)).$$
By the minimality of $R$, we have $R = S \cup\{ u'\}$.
So $u_{1,1} = u'$ and $J[R]$ is a $k$-clique minus the edges $e_1 = u_{1,1}u_{2,1}$ and $e_2 = u_{1,1}u_{2,2}$.
But then for any possible choice of $u_{1,2}$, there exists a $(k-1)$-coloring $\phi$ of $J[R]$ such that $\{\phi'(u_{1,1}), \phi'(u_{1,2})\} \neq \{\phi'(u_{2,1}), \phi'(u_{2,2})\}$.
This contradiction to~\eqref{june246} finishes Case 1.
	
\medskip
In all subsequent cases, we will use Lemma~\ref{lem1 - k} in order to  construct either a $(k-1)$-coloring of $J$ or a $(k-1)$-coloring of $J[R]$
fitting into Case 1. 
Lemma~\ref{lem1 - k} uses variables $i',z,s,w$.  We will use several different values for $i'$, but $z,s,w$ are always defined by the discussion leading to~\eqref{edges to add}, and by Claim~\ref{3-connected}, $s \geq 3$. 

{\bf Case 2:}  
$2i\geq z=|E(R,Q)|$.
By (\ref{se1}) and~\eqref{edges to add}, in order to have $2i\geq |E(R,Q)|$, we need $i=\frac{k}{2}$, $k\geq 6$, and $|E(R,Q)| = k$.
For $k\geq 6$, we know that $y_k=k^2-5k+2$.
By Lemma~\ref{lem1 - k} with $i'=i-1$,   we can add to $J[R_*]$ a set $E_1$ of at most $i-1$ edges
such that (\ref{I01}) holds with $i' =i-1$.
Let $H_1 = J[R] \cup E_1$.
By (\ref{j14'}), 
$$\rho_{k,H_1}(R')>y_k + 2i(k-1) - (i-1)2(k-1) = k(k-3)$$ for every $R'\subseteq R$ with $|R'|\geq 2$ (as this inequality is also true for $K_{k-1}$).
So, by Corollary~\ref{k(k-3)}, $H_1$ has a proper $(k-1)$-coloring $\phi$. 
Moreover, $\phi$ is also a proper $(k-1)$-coloring of $J[R]$.
Since Case 1 does not hold, $\phi$ has a color class
$C$ that satisfies neither~\eqref{low weight color classes}  nor~\eqref{low weight color classes - case 4}.
This means that $\sum_{u \in R_* - C} w(u) = i - 1$ (the lower bound is from (\ref{I01})) and $  \sum_{u \in C} w(u) \geq k-1. $
But then $$|E(R,Q)|\geq k-1+i-1\geq k-1+\frac{k}{2}-1=\frac{3k}{2}-2.$$ Since $k\geq 6$, this contradicts $|E(R,Q)| = k$.
This concludes Case 2.

If Case 2 does not hold, then $z>2i$ and, since $s = |R_*| \geq 3$ by Claim~\ref{3-connected}, 
the ``moreover'' part of Lemma~\ref{lem1 - k} holds when applied with $i'=i$ (recall the calculation right after (\ref{j14'}) showing that $i\geq 1$). The cases below analyze the possibilities in Lemma~\ref{lem1 - k}.

\vspace{3mm}
{\bf Case 3:} The set  $\{w(u_1),\ldots,w(u_s)\}$ is  $i$-special: $s=i+1$ and $w(u_2)=\ldots=w(u_s)=1$.
This means that many (exactly $z - i \geq i$) edges connect $u_1$ with $Q$ and each of the vertices $u_2,\ldots,u_{i+1}$
is connected to $Q$ by exactly one edge.
 For $j=2,\ldots,i+1$, let $q_j$ be the vertex in $Q$ such that $u_jq_j\in E(J)$.
  Let $E_0=\{u_1u_j\,:\, 2\leq j\leq i\}$ and $H_0=J[R]\cup E_0$. 
Since $|R|<n$, $H_0$ is smaller than $J$.
Since $|E_0|=i-1$,   by (\ref{j14'}), 
$$\rho_{k,H_0}(R')>y_k + 2i(k-1) - (i-1)2(k-1) \geq k(k-3)$$ for every $R'\subseteq R$ with $|R'|\geq 2$ (as such an inequality is also true for $K_{k-1}$).
So, by the second part of Corollary \ref{k(k-3)}, $H_0$ has a proper $(k-1)$-coloring $\phi$.
By construction, $\phi$ is a proper $(k-1)$-coloring of $J[R]$ that satisfies $\phi(u_j) \neq \phi(u_1)$ for each $2 \leq j \leq i$.
If $\phi(u_{i+1})\neq \phi(u_1)$, then for every monochromatic subset $M$ of $R_*$ in $J\cup E_0$ with $|M|\geq 2$, (\ref{I01}) holds.
But this coloring satisfies (\ref{low weight color classes}), which contradicts that Case 1 does not apply, so suppose $\phi(u_{i+1})= \phi(u_1)$.

Let $J_0$ be obtained from $J[V(J)-(R-u_1)]$ by adding edge $u_1q_{i+1}$.
By Claim~\ref{very small potential}, $J_0$ has a $(k-1)$-coloring $\phi'$.
 Since $i\leq \frac{k}{2}$ by~\eqref{se1}, we can rename
the colors in $\phi'$ so that $\phi'(u_1)=\phi(u_1)=\phi(u_{i+1})$
and $\phi(\{u_2,\ldots,u_i\})\cap \phi'(\{q_2,\ldots,q_i\})=\emptyset$. Then $\phi\cup \phi'$ is
a proper $(k-1)$-coloring of $J$, a contradiction.

\vspace{3mm}
{\bf Case 4:}  The set of weights $\{w(u_1),\ldots,w(u_s)\}$ is not $i$-special and $2i< z$, so that Part (i) or (ii) of the ``moreover'' part of Lemma~\ref{lem1 - k} holds when applied with $i'=i$.
If Part (i)  holds, then we  take this set
$E_0$ of at most $i-1$ edges and let $H_0=J[R]\cup E_0$. In this case by (\ref{j14'}),
$\rho_{k,H_0}(R')>y_k+2k-2\geq k(k-3)$ for every $R'\subseteq R$ with $|R'|\geq 2$ (as such an inequality  also holds for $K_{k-1}$).
So, by the second part of Corollary \ref{k(k-3)}, $H_0$ has a $(k-1)$-coloring $\phi$, satisfying (\ref{low weight color classes}) of Case 1.

 Suppose now that Part (ii) holds:
{\em there is a hypergraph $H$ with at most $i-1$ graph edges and a $3$-edge $e_0=\{u,v,w\}$ such that
$d_{H}(u_j)\leq w(u_j)$ for all $j=1,\ldots,s$ and (\ref{I01}) holds.} Let $H_1$ be obtained from
$J[R]$ by adding the set of edges $E(H)-e_0$ and edge $uv$. 
Since $|R|<n$, $H_1$ is smaller than $J$.
A proper $(k-1)$-coloring of $H_1$ would satisfy (\ref{low weight color classes}) of Case 1, so $\chi(H_1) \geq k$.
Then $H_1$ has a $k$-critical subgraph $H'_1$. 

Let $R'=V(H'_1)$. 
Because Part (i) of the ``moreover'' does not hold,  $\{u,v\} \subseteq R'$ and $uv \notin J[R']$. So $J[R']$ is not a clique.
If $H'_1$ is not a $k$-Ore graph, then by the minimality of $J$, $\rho_{k,H_1}(R')\leq y_k$ and so
$\rho_{k,J}(R')\leq y_k+2i(k-1)$, contradicting in~\eqref{j14'} the minimality of $i$ (as $y_k+2i(k-1) \leq 2(k-1)(k-2)$ and $J[R']$ is not a clique).
Thus, $H'_1$ is  a $k$-Ore graph and
 $\rho_{k,H_1}(R') = k(k-3)\leq y_k+2k-2$.  Then
$$\rho_{k,J}(R')\leq \rho_{k,H_1}(R')+2i(k-1)\leq y_k+2(i+1)(k-1),$$ and by the minimality of $R$, 
$\, R'=R$.
Furthermore,   if $H'_1\neq H_1$, then it has the same vertex set as $H_1$ and at least one
fewer edge, in which case,
$$\rho_{k,J}(R')\leq \rho_{k,H'_1}(R')+2i(k-1)\leq  \rho_{k,H_1}(R')+2(i-1)(k-1)\leq k(k-3)
+2(i-1)(k-1)\leq y_k+2i(k-1),$$
a contradiction to (\ref{j14'}). So, $H_1'=H_1$, $H_1$ is a $k$-Ore graph and so $\rho_{k,J}(R)=k(k-3)+2i(k-1)$.

Recall that $\rho_{k,J}(R) \leq 2(k-1)(k-2)$ by the last inequality of (\ref{j14'1}), so 
$$i \leq \frac{2(k-1)(k-2) - k(k-3)}{2(k-1)} = \frac{1}{2}(k-2 + \frac{2}{k-1}).$$
Because $i$ is an integer, this inequality is strict, and so $\rho_{k,J}(R) < 2(k-1)(k-2)$. 
Therefore $\rho_{k,J}(W) > \rho_{k,J}(R)$ if $J[W] = K_{k-1}$.
By the minimality of $i$ and $R$, any $W \subset R$ such that $|W| \geq 2$ and $J[W] \neq K_{k-1}$ satisfies $\rho_{k,J}(W) > \rho_{k,J}(R)$.
Graph $H_1 - uv$ is $J[R]$ plus $i-1$ edges, so for any $W \subset V(H_1')$ with $|W| \geq 2$ we have 
$$\rho_{k,H_1-uv}(W) = \rho_{k,J}(W) -2(i-1)(k-1) > k(k-3) + 2(k-1) = (k+1)(k-2).$$
Thus by Lemma~\ref{extr}, $H_1-uv$ has a $(k-1)$-coloring $\phi$ with $\phi(w)\neq \phi(u)$. 
This is a $(k-1)$-coloring of $J[R]$ that satisfies \eqref{low weight color classes} of Case 1.
\qed

Recall that a standard set has potential $(k+1)(k-2)$.
Because $(k+1)(k-2) < 2(k-1)(k-2)$ when $k \geq 4$, Lemma~\ref{small potential} implies: 

\begin{cor}\label{no standard set} For each $R \subsetneq V(J)$ with $2\leq |R|\leq n-2$,  $\rho_k(R) > (k+1)(k-2)$.
In particular, $J$ does not contain  a standard set of size at most $n-2$.
\end{cor}

\begin{claim} \label{small clusters}
If $v$ is not in a $(k-1)$-clique $X$, then $|N(v) \cap X| \leq \frac{k-1}{2}$.
Furthermore, if $T$ is a cluster in a $(k-1)$-clique $X$, then $|T| \leq \frac{k-1}{2}$.
\end{claim}
{\bf Proof.}
If $|N(v) \cap X| \geq \left\lceil k/2\right\rceil$, then $\rho_k(X + v) \leq 2(k-2)(k-1) - 2$.
Since $n\geq k+2$, this  contradicts Lemma~\ref{small potential}. This proves the first part.

Suppose now that $T$ is a cluster in a $(k-1)$-clique $X$. Since $|X|=k-1$ and $d(w)=k-1$  for every $w\in T$, 
each such $w$ has the unique neighbor  $v(w)$ outside of $X$. But by the definition of a cluster, $v(w)$ is
the same, say $v$, for all $w\in T$.
This means that $T \subseteq X \cap N(v)$, so $|N(v) \cap X| \geq |T|$. Thus the second part follows from the first.
\qed

\begin{claim} \label{big neighbors (b)}
Suppose $T$ is a cluster in $J$, $t=|T| \geq 2$, and  $N(T)\cup T$ contains a $(k-1)$-clique $X$.
Then  $d_J(v)\geq k-1+t$ for every $v\in X-T$.
\end{claim}
{\bf Proof.} 

By the definition of a cluster,   $|T \cup N(T)| = k$. So, since $|T|\geq 2$ and $|X|=k-1$,
$T\cap X\neq \emptyset$. So by Claim~\ref{unique cluster}(ii), $T\subseteq X$,
and by Claim~\ref{unique cluster}(i), $X$ contains only one nonempty cluster, namely, $T$.

Suppose  $v\in X-T$ and $d(v) \leq k-2 + t$. By the previous paragraph,
$v$ is not in a cluster and thus by Fact \ref{clusters partition}, $d(v) \geq k$.
 By  Claim~\ref{minus edge}, $T$ is contained in at most one $(k-1)$-clique (which is $X$), and so
\begin{equation}\label{j15}
\mbox{ $N(T)\cup T-v$ does not contain $K_{k-1}$.}
\end{equation}
Because $T$ and $v$ are parts of the same clique, $|N(v) - T| = d(v) - |T|$, and  this is at most $k-2$, 
since $d(v) \leq k-2 + t$.
Let $u \in T$ and $J' = J - v + u'$, where $u'$ is a new vertex that satisfies $N[u'] = N[u]$.
Suppose $J'$ has a $(k-1)$-coloring $\phi':V(J') \rightarrow C =\{c_1, \dots c_{k-1}\}$.
Then we find a $(k-1)$-coloring $\phi$ of $J$ as follows:
set $\phi|_{V(J) - T - v} = \phi'|_{V(J') - T - u'}$, $\phi(v) \in C - \phi'(N(v) - T)$,
and then color $T$ using colors in $\phi'(T \cup u') - \phi(v)$.
This is a contradiction, so there is no $(k-1)$-coloring of $J'$.
Thus $J'$ contains a $k$-critical subgraph $J''$.
Let $W=V(J'')$. 
By Corollary \ref{k(k-3)}, $\rho_{k,J'}(W) \leq k(k-3)$.

By the criticality of $J$, graph $J''$ is not a subgraph of $J$.
So $u' \in W$ and $u'$ has at least $k-1$ neighbors in $W$.
By symmetry, we have $T \subset W$.
But then
$$ \rho_{k,J}(W-u') \leq k(k-3) - (k-2)(k+1) + 2(k-1)(k-1) = 2(k-2)(k-1).$$
By Lemma \ref{small potential}, this implies that either $J[W - u']$ is a $K_{k-1}$ or $W-u'=V(J)-v$. 
If the former holds, then because $J[W - u']$ is a complete graph and $T \subset W-u'$ 
we have  $N(T) \cup T \supset J[W - u'] \cong K_{k-1}$,  and because $v \notin W$ this is a contradiction to (\ref{j15}). 
If the latter holds, then we have a contradiction to Fact \ref{aug3}, since $d(v)\geq k$.
 \qed

\begin{claim} \label{adjacent k-1}
Suppose $xy \in E(J)$, $N[x] \neq N[y]$, $x$ is in a cluster of size $s$, $y$ is in a cluster of size $t$, and $s \geq t$.
Then $x$ is in a $(k-1)$-clique.
Furthermore, $t = 1$.
\end{claim}
{\bf Proof.}
Assume that $x$ is not in a $(k-1)$-clique.
Let $J' = J - y + x'$, where $x'$ is a new vertex with $N_{J'}(x')=N_J[x]-y$.
By the definition of a cluster,  $d(x) = d(y) = k-1$.
Both $J'$ and $J$ have the same number of vertices and the same number of edges (because $xy \in E(J)$, vertex
 $x$ lost an edge to $y$ and gained an edge to $x'$). 
If distinct $v,w\in V(J)-y$, then they have the same closed neighborhood in $J'$ if and only if they have the same closed neighborhood in $J$.
By the definition of a cluster, there are exactly $t-1$  vertices in $V(J)- y$ that have the same closed 
neighborhood as $y$.
By construction, there are exactly $s$  vertices in $V(J')-x'$  that have the same closed neighborhood in $J'$ as $x'$.
So since $s\geq t$,  by Rule (S3), $J'$ is smaller than $J$.

If $J'$ has a $(k-1)$-coloring $\phi':V(J') \rightarrow C = \{c_1, c_2, \dots c_{k-1}\}$, 
then we extend it to a proper $(k-1)$-coloring $\phi$ of $J$ as follows: define $\phi|_{V(J)-x-y} = \phi'|_{V(J')-x-x'}$, 
then  choose $\phi(y) \in C - (\phi'(N(y) - x))$, and $\phi(x) \in \{\phi'(x), \phi'(x')\} - \{\phi(y)\}$.

So, $\chi(J')\geq k$ and $J'$ contains a $k$-critical subgraph $J''$.
Let $W=V(J'')$. 
By the criticality of $J$ and because $y \notin J''$, we have  $J'' \neq J$, and $J''$ is not a subgraph of $J$.
This yields $x' \in W$.
By symmetry, $x \in W$. By Corollary~\ref{k(k-3)}, $\rho_{k,J'}(W)\leq k(k-3)$.
Because of this and the fact that $d(x') = k-1$, 
\begin{equation}\label{adjacent k-1 cluster}
\rho_{k,J}(W - x') \leq k(k-3) - \rho_{k,J'}(\{x'\}) + 2(k-1)d(x') = 2(k-2)(k-1).
\end{equation}
By assumption, $x$ is not in a $(k-1)$-clique, so Lemma \ref{small potential} implies that $|W-x'| > n-2$.
Since $y \notin W$, this yields
 $|W-x'| = n-1$, which implies $V(J') = V(J'')$ and  $W-x' = V(J)-y$.
By Corollary \ref{k(k-3)}, $\rho_{k,J''}(W) \leq k(k-3)$.
Moreover, because $J''$ is smaller than or equal to $J'$ which is smaller than $J$,  the minimality of $J$ implies 
that if $J''$ is not $k$-Ore then $\rho_{k,J''}(W) \leq y_k$. In this case we can replace the term $k(k-3)$ in
(\ref{adjacent k-1 cluster}) with $y_k$ (as $J''$ is a subgraph of $J'$ and so $\rho_{k,J'}(W) \leq \rho_{k,J''}(W)$) and get $\rho_{k,J}(W - x')\leq y_k+k^2-3k+4$ contradicting Fact \ref{aug3}.
If $J'' \neq J'$ then $\rho_{k,J''}(W) - 2(k-1) \geq \rho_{k,J'}(W)$, and instead of (\ref{adjacent k-1 cluster})
we obtain $\rho_{k,J}(W - x') \leq  2(k-2)(k-1)-2(k-1)$, again contradicting Fact \ref{aug3}.
So  $J'' = J'$ and $J''$ is $k$-Ore, thus $J'$ is a $k$-Ore graph.

Since $n>k$, $J'\neq K_{k}$.
Let the  separating set $\{u,v\}$, vertex subsets $A=A(J',u,v)$ and $B=B(J',u,v)$, and graphs 
$ \widetilde J'(u,v)$ and $\check J'(u,v)$  be as in Fact~\ref{f2}.
By Corollary~\ref{no standard set}, the standard set $A$  is not contained in $V(J)$. Hence $x' \in A$.
Therefore $x' \notin V(\check J'(u,v))-u*v$. Thus for every  $W\subseteq V(\check J')-u*v$ with $|W|\geq 2$,
by Corollary \ref{no standard set},
we have $\rho_{k,\check J'(u,v)}(W) = \rho_{k,J}(W)>(k+1)(k-2)$.
Then by Claim~\ref{special vertex}, there exists an 
\begin{equation}\label{j11}\mbox{\em
$S \subseteq V(\check J'(u,v))-u*v$ such that 
$\check J'(u,v)[S] \cong K_{k-1}$, and $d_{\check J'(u,v)}(w) = k-1$ for all $w \in S$.}
\end{equation}
By Claim~\ref{minus edge}, vertex $y$ in $J$ is adjacent to at most $k-3$ vertices in $S$.
Because $V(J)-y\subset V(J')$,~\eqref{j11} yields that vertices in $S-N(y)$ have degree $k-1$ in $J$.
So by Fact \ref{clusters partition}, $S$ contains a $V(J)$-cluster $T$, which by Claim~\ref{unique cluster} contains all vertices in $S$ of degree $k-1$ and $|T|\geq|S-N(y)|\geq 2$.
Then by Claim~\ref{big neighbors (b)}, the degree of each vertex in $S-T$ in $J$ is at least $k+1$.
This is impossible, since by Fact \ref{f1.1}.5 each of them has in $J$ at most one extra neighbor (and it is $y$)
in comparison with $\check J'(u,v)$, where their degree was only $k-1$.
This proves the first part:  $x$ is in a $({k-1})$-clique, say $X$.

Let $T_y$ be the cluster containing $y$.
By the definition of a cluster, every vertex in $T_y$ has the same neighbors as $y$, and so $T_y \subseteq N(x)$.
The $(k-1)$-clique $X$ containing $x$ is a part of  $N[x]$.
By Claim~\ref{unique cluster}(ii) if $T_y \cap X \neq \emptyset$, then $T_y \subseteq X$, and by Claim~\ref{unique cluster}(i) $y \in T_y - X$; so $T_y \cap X = \emptyset$.
Therefore $|T_y| \leq |N(x) - X| = d(x) - (k-2) = 1$. This proves the second part.
\qed

\begin{claim} \label{big neighbors (a)}
Suppose $T$ is a cluster in $J$, $t=|T| \geq 2$, and $N(T)\cup T$ does not contain $K_{k-1}$.
Then  $d_J(v)\geq k-1+t$ for every $v\in N(T)-T$.
\end{claim}
{\bf Proof.}
By Claim~\ref{adjacent k-1},  $d(v)\geq k$.
Now the proof follows exactly as the proof to Claim \ref{big neighbors (b)}.
\qed

\section{Proof of Theorem \ref{pot theorem}}
Now we are ready to prove the theorem.
Recall that in Section $3$ we defined a relation ``smaller,'' and 
 $J$  is a ``smallest'' counterexample to Theorem \ref{pot theorem}: it is a $k$-critical 
	graph with
	 $\rho_k(V(J)) > y_k$ and is not $k$-Ore.

We will use the following result on $k$-critical graphs which is Corollary~9 in~\cite{KY}.
\begin{lemma}[\cite{KY}]\label{co1} Let $J$ be a $k$-critical graph. Let
 disjoint vertex subsets  $A$ and $B$ be such that\\
 (a)  at least one of $A$ an $B$ is independent;\\
(b) $d(a)=k-1$ for every $a\in A$;\\
(c) $d(b)=k$ for every $b\in B$;\\
(d) $|A|+|B|\geq 3$.\\
 Then (i) $e(J(A,B))\leq 2(|A|+|B|)-4$ and
(ii) $e(J(A,B))\leq |A|+3|B|-3$.
\end{lemma}


\subsection{Case $k = 4$}
In this subsection we prove the theorem for $k=4$.
Specifically, we will prove that $|E(J)| \geq \frac{5}3 |V(J)|$, which will imply that $\rho_{4,J}(V(J)) \leq 
 0\leq \ y_4 = 2$.

\begin{claim} \label{two 3 etc}
Each vertex with degree $3$ has at most one neighbor with degree $3$.
\end{claim}
{\bf Proof.}
Let $x$ be such that $N(x) = \{a,b,c\}$ and $d(a) = 3$.
Then $x$ and $a$ are each in a cluster.
Because no cluster is larger than $k-3=1$ by Claim \ref{unique cluster}, $a$ and $x$ are in different clusters.
Then by Claim \ref{adjacent k-1}, $J[\{x,b,c\}]$ is a $K_3$.
So by Claim~\ref{unique cluster}, $d(b), d(c) \geq 4$.
\qed

We  now use discharging to show that $|E(J)| \geq \frac{5}3 n$.
Each vertex begins with charge equal to its degree.
If $d(v) \geq 4$, then $v$ gives charge $\frac16$ to each neighbor.
Note that $v$ will be left with charge at least $\frac56 d(v) \geq \frac{10}{3}$.
By Claim \ref{two 3 etc}, each vertex of degree $3$ will end with charge at least $3 + \frac26 = \frac{10}3$.
Therefore the total charge is at least $\frac{10}3n$, and thus so is the sum of the vertex degrees.
Hence the number of edges is at least $\frac{5}3n$.
\qed

\subsection{Case $k = 5$}
In this subsection we prove the theorem for $k=5$.
Specifically, we will prove that $|E(J)| \geq \frac{9}4 |V(J)|$, which will imply that $\rho_{5,J}(V(J)) \leq 0<y_5 = 4$.

\begin{claim} \label{4 no neighbors}
If $k = 5$, then each cluster has only one vertex.
\end{claim}
{\bf Proof.} Suppose the claim does not hold. 
By Claim \ref{unique cluster}, every cluster has size at most $k-3 = 2$, so assume that $\{x,y\}$ is a cluster: $N[x] = N[y]$ and $d(x) = d(y) = 4$. 
Let $N(x)=\{y,a,b,c\}$. 
By assumption $J$ is not $5$-Ore and therefore $J\neq K_5$ (and since $J$ is critical, it does not contain a $K_5$).
By Claim \ref{very small potential}, $J$ does not contain a subgraph isomorphic to $K_5-e$.
Therefore any five vertices in $J$ induce at most ${5 \choose 2} - 2$ edges, and thus $|E(J[\{a,b,c\}])|\leq 1$.
By Claims \ref{big neighbors (b)} and \ref{big neighbors (a)}, we can rename the vertices in $\{a,b,c\}$ so that $ab,ac\notin E(J)$ and $d(c)\geq 6$.

We obtain $J'$ from $J$ by deleting $x$ and $y$ and gluing $a$ with $b$.
If $J'$ is $4$-colorable, then so is $J$.
This is because  any $4$-coloring of $J'$ will have at most $2$ colors on $N[x] - \{x,y\}$
and therefore could be extended greedily to $x$ and $y$.

So $J'$ contains a $k$-critical subgraph $J''$.
Let $W'=V(J'')$. 
Then by Corollary \ref{k(k-3)}, $\rho_{5,J'}(W') \leq \rho_{5,J''}(W') \leq 10$.
Furthermore, because $J''$ is smaller than $J$, 
\begin{equation}\label{J1}
\mbox{if $J''$ is not $k$-Ore, then  $\rho_{5,J''}(W') \leq 4$.}
\end{equation}

Because $J$ is critical and $x,y \notin J'' \subseteq J'$, graph  $J''$ is not a subgraph of $J$.
This implies that $a*b \in V(J'')$.
Let $W = W' - a*b + a + b + x + y$.
By the definition of  potential, 
\begin{equation}\label{k=5 cluster}
\rho_{5,J}(W) = \rho_{5,J''}(W') + 18\left(|W|-|W'|\right) - 8\left(|E(J[W])|-|E(J''[W'])|\right).
\end{equation}
If $c\notin W'$, then by construction, $W$ has $3$ more vertices and induces at least $5$ more edges than $W'$.
If $c\in W'$, then $W$ has $3$ more vertices and induces at least $7$ more edges than $W'$.
Suppose first that $c \notin W'$, so that~\eqref{k=5 cluster} becomes $\rho_{5,J}(W) \leq 10 + 54 - 40 = 24$.
Because $ab \notin E(J)$, $J[W]$ is not a $K_4$.
So by Lemma \ref{small potential}, $|W|\geq n-1$.
Therefore $W = V(J) - c$ and $\rho_{5,J}(W) \leq 24$, but this contradicts Fact \ref{aug3} because $d(c) \geq 6 = k+1$.

So now we assume that $c \in W'$, which means that~\eqref{k=5 cluster} becomes $\rho_{5,J}(W) \leq 10 + 54 - 56 \leq 8$.
By Lemma~\ref{very small}, $W = V(J)$, which then implies  $V(J'') = V(J')$.
Furthermore, if $J''$ is not $k$-Ore, then by~\eqref{J1},  $\rho_{5,J''}(W')\leq 4$ and the right hand side of~\eqref{k=5 cluster} has an extra $-6$.
If $J''$ is a proper subgraph of $J'$, then $|E(J[W])|-|E(J''[W'])|$ is at least $8$ instead of the previously applied bound of $7$.
In either case our bound becomes stronger by at least $6$, and becomes $\rho_{5,J}(V(J)) = \rho_{5,J}(W) \leq 8 - 6 = 2 < y_k$, a contradiction to the choice of $J$.
So $J''$ is $k$-Ore and $J''=J'$.

This also implies that $N(a) \cap N(b) = \{x,y\}$, because otherwise, since $J'' = J'$,   
we would have gained an extra edge when we undo the merge of $a$ and $b$ into $a*b$, 
which would increase $|E(J[W])|-|E(J''[W'])|$ by $1$ in~\eqref{k=5 cluster} 
yielding  the same contradiction.

Since $d(c) \geq 6$,  $J'$ cannot be $K_5$.
Let the separating set $\{u,v\}$, vertex subsets $A=A(J',u,v)$ and $B=B(J',u,v)$, and graphs $ \widetilde J'(u,v)$ and $\check J'(u,v)$  be as in Fact~\ref{f2}.
By Fact \ref{f1.1}.3, $\rho_{k,J'}(A) = (k+1)(k-2)$. 
By Corollary \ref{no standard set}, proper subgraphs of $J$ have strictly larger potential than $(k+1)(k-2)$ and so $a*b \in A$.
Therefore $J'[B - u- v] \subset J$ and so $V(\check J') - V(J) = u*v$.
Hence by  Corollary~\ref{no standard set}, $\rho_{k,\check J'(u,v)}(U) = \rho_{k,J}(U) >(k+1)(k-2)$ 
for every  $U\subseteq V(\check J')-u*v$ with $|U|\geq 2$.
Then by Claim~\ref{special vertex}, there exists an $S \subseteq V(\check J'(u,v))-u*v$ such that $\check J'(u,v)[S] \cong K_{k-1}$, 
and $d_{\check J'(u,v)}(w) = k-1 = 4$ for all $w \in S$.

We claim that 
\begin{equation}\label{j142}
 \mbox{\em  each vertex in $S-c$ has degree $k-1$ in $J$.}
\end{equation}
Indeed, how it is possible that a vertex  $w\in S$  has a larger degree in $J$ than in $\check J'$?
By Fact \ref{f1.1}.7, $d_{\check J'}(w)=d_{J'}(w)$.  
Because $N(a) \cap N(b) = \{x,y\}$, the only vertices whose degree in $J$ could be greater than in  $J'$  are $a,b,c$.
But we already showed that $a*b \notin V(\check J')-u*v$ and $S \subseteq V(\check J'(u,v))-u*v$, so  $S \cap \{a,b,c\} \subseteq \{c\}$.
This proves~\eqref{j142}.

By Fact \ref{clusters partition} and~\eqref{j142},
 every vertex in $S-c$ is in a cluster.
Because $S$ is a $(k-1)$-clique, by Claim \ref{unique cluster} there is only one cluster in $S$, so the claim
implies that $T = S-c$ is a cluster and $|T|= 3$, which contradicts Claim \ref{unique cluster} 
that each cluster in $J$ has size at most $k - 3 = 2$.
\qed

\begin{claim}\label{c23}
Each $K_4$-subgraph of $J$ contains at most one vertex with degree $4$.
Furthermore, if $d(x) = d(y) = 4$ and $xy \in E(J)$, then each of $x$ and $y$ is in a $K_4$.
\end{claim}
{\bf Proof.}
Each vertex of degree $4$ is in a cluster by definition, and by Claim \ref{unique cluster}, each $K_4$ contains only one cluster.
The first statement of our claim then follows from Claim \ref{4 no neighbors}  
and the second --- from Claim~\ref{adjacent k-1}.
\qed

\begin{defn}
Let $H \subseteq V(J)$  be the set of vertices of degree $5$ not in a $K_4$,
and $L \subseteq V(J)$  be the set of vertices of degree $4$ not in a $K_4$.
Set $\ell=|L|$, $h=|H|$ and $e_0 = |E(L, H)|$.
\end{defn}

\begin{claim} \label{HL charge}
$e_0 \leq  3h + \ell$.
\end{claim}
{\bf Proof.}
This is trivial if $h+\ell\leq 2$.
By Claim \ref{c23}, $L$ is  independent.
So the claim follows by  Lemma~\ref{co1}(ii) with $A = L$ and $B = H$.
\qed

We will now use discharging to show that $|E(J)| \geq \frac{9}4 n$, which will finish the proof to the case $k=5$.
Let every vertex  $v\in V(J)$ have initial charge $d(v)$.
The discharging has one rule:

{\bf Rule R1:} Each vertex in $V(J)-H$ with degree at least $5$ gives charge $1/6$ to each neighbor.

We will show that the charge of each vertex in $V(J) - H - L$ is at least $4.5$, and then show that the average charge of the vertices in
$H \cup L$ is at least $4.5$.

\begin{claim} \label{5 discharging 1}
After discharging, each vertex in $V(J) - H - L$ has charge at least $4.5$.
\end{claim}
{\bf Proof.}
Let $v \in V(J) - H - L$.
If $d(v) = 4$ and $v \notin L$, then $v$ is in a $K_4$ and by Claim \ref{c23} $v$ receives charge $1/6$ from at least $3$ neighbors and gives no charge.
If $d(v) = 5$ and $v \notin H$, then $v$ is in a $K_4$ and by Claim \ref{c23} $N(v)$ contains at least $2$ vertices with degree at least $5$.
Therefore $v$ gives charge $1/6$ to $5$ neighbors, but receives charge $1/6$ from at least $2$ neighbors.
If $d(v) \geq 6$, then $v$ is left with charge at least $5d(v)/6 \geq 4.5$.
\qed

\begin{claim} \label{5 discharging 2}
After discharging, the sum of the charges on the vertices in $H \cup L$ is at least $4.5|H \cup L|$.
\end{claim}
{\bf Proof.}
By Claim~\ref{c23}, if $v \in L$ then every vertex in $N(v)$ has degree at least $5$.
By Rule R1,  the vertices in $L$ receive from outside of $H\cup L$ the charge at least $\frac{1}{6}(4\ell-|E(H,L)|)$.
By Claim \ref{HL charge}, $|E(H,L)| \leq  3h + \ell $.
So, the total charge on $H\cup L$ is at least
$$5h+4\ell+\frac{1}{6}(4\ell-(3h+\ell))=4.5(h+\ell),$$
as claimed.
\qed

Combining Claims \ref{5 discharging 1} and \ref{5 discharging 2}, the total charge is at least $\frac{9}2n$.
Thus the  degree sum of $V(J)$  is at least $\frac{9}2n$, and so $|E(J)| \geq \frac94|V(J)|$. \qed

\subsection{Case $k \geq 6$}
In this subsection we prove Theorem \ref{pot theorem} for $k\geq6$.
We will prove that $|E(J)| \geq \frac{(k+1)(k-2)}{2(k-1)} |V(J)|$, which will imply that $\rho_{k,J}(V(J)) \leq 0 \leq y_k = k^2 - 5k + 2$.
This proof will involve several claims. 

\begin{claim} \label{cliques have k+1}
Suppose $k \geq 6$, $X$ is a $(k-1)$-clique, and $v\in X$ has degree $k-1$.
Then $X$ contains at least $(k-1)/2$ vertices with degree at least $k+1$.
\end{claim}
{\bf Proof.}
Let $\{u\} = N(v) - X$.
Assume that $X$ contains at least $k/2$ vertices with degree at most $k$.
By Claim \ref{small clusters}, $|N(u) \cap X| < k/2$, so there exists a $w \in X$ such that $uw \notin E(J)$ and $d(w) \leq k$.
By Claim \ref{one cluster}, $d(w) = k$, so assume $N(w) - X = \{a,b\}$.
Let $J'$  be obtained from $J-v$ by adding edges $ua$ and $ub$ if they do not already exist.

Suppose $J'$ has a $(k-1)$-coloring $f$.
If $f(u)$ is not used on $X-w-v$, then we  recolor $w$ with $f(u)$.
So, $v$ will have at least two neighbors of color $f(u)$, and we can extend the $(k-1)$-coloring to $v$.

Thus $J'$ is not $(k-1)$-colorable and so  contains a $k$-critical subgraph $J''$.
Let $W=V(J'')$.
By Corollary \ref{k(k-3)}, $\rho_{k,J'}(W) \leq k(k-3)$ and so 
$$\rho_{k,J}(W) \leq k(k-3) + 2(k-1)(2)=k^2+k-4 < 2(k-2)(k-1).$$
If $W \neq V(J')$ then this contradicts Lemma \ref{small potential}, since in this case $|W|\leq |V(J')|-1\leq n-2$.
So,  $W = V(J')$.

If $J''$ is not a $k$-Ore graph, then by the minimality of $J$, $\rho_{k,J''}(W)\leq y_k$, 
and since  adding edges only reduces potential, we have $\rho_{k,J'}(W) \leq \rho_{k,J''}(W)$, and so
$$\rho_{k,J}(V(J)) \leq \rho_{k,J'}(W) + (k-2)(k+1)(1) - 2(k-1)(k-3) < y_k$$
when $k \geq 6$. This contradicts our choice of $J$, and so $J''$ is  $k$-Ore.

Suppose $J'' \neq J'$. 
Since $W = V(J')=V(J'')$, it follows that $\rho_{k,J'}(W) \leq \rho_{k,J''}(W) - 2(k-1)$, which leads to the same contradiction 
because  by Fact \ref{f1}.8, $\rho_{k,J''}(W) - 2(k-1) \leq k(k-3) - 2(k-1) \leq y_k$.
So, our case is that $J''$ is a $k$-Ore graph and $J'' = J'$, i.e.,  $J'$ is a $k$-Ore graph.
Therefore $\rho_{k,J''}(W) = k(k-3)$ by Fact \ref{f1}.10.

Since $J'-ua-ub$ is a subgraph of $J$, by Corollary \ref{no standard set}, $\rho_{k,J'}(U) = \rho_{k,J}(U) >(k+1)(k-2)$ for every  $U\subseteq V(J')-u$ with $|U|\geq 2$.
Then by Claim~\ref{special vertex}, there exists an
\begin{equation}\label{1210}
\mbox{ \em $S \subseteq V(J')-u$ such that $J'[S] \cong K_{k-1}$, and  $d_{ J'}(v) = k-1$ for all $v \in S$.}
\end{equation}
Note that for every $z\in S-a-b-N_J(v)$, we have $d_J(z)=d_{J'}(z)=k-1$.

Suppose $z\in S \cap N_J(v)$. Then since $N_J(v) = \{u\} \cup X$ and  $u \notin S$, we have $z \in X\cap S$.
So $z \notin \{a,b\}$ and therefore $N_{J'}(z) = N_J(z) \cup \{v\}$, and thus $d_J(z) \geq |X \cup S| - 1$.
Because $|S| = k-1$, there exists a $z' \in S - X$.
By construction of $J'$ and because $S$ is a clique, we have that $X \cap S \subseteq N_J(z') \cap X$.
By Claim~\ref {small clusters},  $|X\cap S|\leq|N_J(z') \cap X|\leq (k-1)/2$. 
So for $k\geq 6$, 
$$d_{J'}(z) = d_J(z)-1 \geq |X\cup S|-2\geq \lceil\frac{k-1}{2}\rceil+(k-1)-2\geq k,$$
a contradiction to~\eqref{1210}. Thus, $S \cap N(v) = \emptyset$.

By Fact~\ref{clusters partition} and Claim~\ref{unique cluster}, the vertices in $S-a-b$ are part of one cluster in $J$ with size at least $k-1-2$ in the $(k-1)$-clique $S$, which contradicts Claim~\ref {small clusters} since $k-3 > \frac{k-1}2$ for $k \geq 6$.
 \qed

\begin{claim} \label{k=6}
If $k=6$ and a cluster $T$ is contained in a $5$-clique $X$, then $|T| = 1$.
\end{claim}
{\bf Proof.}
By Claim \ref{small clusters}, assume that $T = \{v_1, v_2\}$.
Let $N(v_1) - X = \{y\}$ and $T-X=\{u, u', u''\}$.
By Claim \ref{cliques have k+1}, $d(u),d(u'),d(u'') \geq 7$.
By Claim \ref{small clusters}, $|N(y) \cap X| < k/2$, and by the definition of a cluster, $\{v_1,v_2\}\subset N(y)$. 
 Thus, $N(y) \cap \{u, u', u''\} = \emptyset$.
Obtain $J'$ from $J-T$ by gluing $u$ to $y$.

Suppose that $J'$ has a $5$-coloring.
Then we can extend this coloring to a $5$-coloring of $J$ by greedily assigning colors to $T$, because only $3$ different colors appear on the set $\{u, u', u'', y\}$.
So we may assume that $\chi(J') \geq 6$.
Then $J'$ contains a $6$-critical subgraph $J''$.
Let $W=V(J'')$. Then by Corollary~\ref{k(k-3)}, $\rho_{6,J'}(W) \leq 6(6-3)= 18$.
Since $J''$ is not a subgraph of $J$ because $J$ itself is critical, $u*y \in W$.
Let $t = | \{u', u''\} \cap W|$.

{\bf Case 1:} $t=0$.
Then 
\begin{equation}\label{162}
\rho_{6, J}(W-u*y + y + X) \leq 18 + 28(5) - 10(12) = 38,
\end{equation}
By Lemma \ref{small potential}, $|W-u*y + y + X |\geq n-1$.
The bound~\eqref{162} did not account for the edges in $E(\{u', u''\}, V(J)-X)$, but each of $u',u''$ has at least $3$ neighbors outside of $X$.
Thus we strengthen~\eqref{162} to $\rho_{6, J}(W-u*y + y + X) \leq 38-10 \cdot 4 < 0$.

\smallskip 
Denote $R=W-u*y + y + u + T$.

{\bf Case 2:} $t = 1$. 
Then $\rho_{6, J}(R) \leq 18 + 28(3) - 10(7) =32$.
By Lemma \ref{small potential}, this yields $|R|\geq n-1$, so $R$ is either $V(J)-u'$ or $V(J)-u''$.
But because $d(u'),d(u'') \geq 7 = k+1$, Fact \ref{aug3} says that $\rho_{6,J}(V(J) - u'), \rho_{6,J}(V(J) - u'') > y_k + k^2 + k = 50$, 
 a contradiction.

{\bf Case 3:} $t = 2$. Then $|R|\geq 4$ and
\begin{equation}\label{163}
 \rho_{6, J}(R) =\rho_{k,J''}(W) + 28\left(|R|-|W|\right) - 10\left(|E(J[R])| - |E(J''[W])| \right)  \leq 18 + 28(3) - 10(9) =12.
 \end{equation}
 Since $|R|\geq 4$,
by Lemma \ref{small potential},  $|R|\geq n-1$. 
  By Fact \ref{aug3}, every $(n-1)$-element set has potential more than $y_k + k^2 - 3k + 4 = 26$.
So  $|R|=n$ and  $R=W-u*y + y + u + T = V(J)$.

If $J''$ is not $k$-Ore, then by the minimality of $J$, $\rho_{k,J''}(W) \leq y_k =8$.
 In this case, we can replace in~\eqref{163} the upper bound of $18$ on $\rho_{k,J''}(W)$ with $8$,
which yields $\rho_{6, J}(R) \leq 2<y_6$, a contradiction. Thus  $J''$ is  $k$-Ore.

Similarly, 
if $J'' \neq J'$, then because $V(J')=V(J'')$ we can replace in~\eqref{163} the term $-10(9)$ with $-10(10)$ as there is an extra edge in $|E(J[R])| - |E(J''[W])|$, which again  yields $\rho_{6, J}(R) \leq 2<y_6$.
So $J' = J''$.

Since $J''-u*y$ is a subgraph of $J$, by Corollary \ref{no standard set},  $\rho_{k,J'}(U) = \rho_{k,J}(U) >(k+1)(k-2)$ for every  $U\subseteq V(J')-u*y$ with $|U|\geq 2$.
Then by Claim~\ref{special vertex}, there exists an $S \subseteq V(J')-u*y$ such that $J'[S] \cong K_{5}$, and  $d_{ J'}(v) = 5$ for all $v \in S$.
By construction, $S \subseteq V(J)$, so $J[S] \cong K_5$.
By Fact \ref{clusters partition} each vertex with degree $k-1$ in $J$ is in a cluster, by Claim \ref{unique cluster} the $S$-cluster is unique, 
and by Claim~\ref {small clusters},  $S$ has at most $2$ vertices of degree $k-1$ in $J$.
So  there are at least three vertices $z_1, z_2, z_3 \in S$ such that $d_J(z_i) > d_{J'}(z_i)$ for $1 \leq i \leq 3$.
But the  vertices with larger  degrees  in $J$ than in $J'$ could be of only two types: (a) those in $N_J[v_1] = N_J[v_2] = \{y, v_1, v_2, u, u', u''\}$ and (b) those in $N_J(y) \cap N_J(u)$.
Because $J'' - u*y$ does not contain $T = \{v_1, v_2\}$ (which was deleted), and $u*y \notin S$, we have at most two vertices of type (a).
Then we  have  a vertex of type (b), but in this case we have an extra edge whose potential was not accounted in~\eqref{163}.
So in this case we again instead of~\eqref{163} get $\rho_{6, J}(R) \leq 12-10<y_6$,
a contradiction.
\qed

\begin{defn} \label{L,H,e0}
We partition $V(J)$ into four classes: $L_0$, $L_1$, $H_0$, and $H_1$.
Let $H_0$ be the set of vertices with degree $k$, $H_1$ be the set of vertices with degree at least $k+1$, and $H = H_0 \cup H_1$.
Let
$$L = \{u \in V(J): d(u) = k-1\},
$$
$$L_0 = \{u \in L : N(u) \subseteq H\},
$$
and
$$ L_1 = L - L_0.
$$
Set $\ell=|L_0|$, $h=|H_0|$ and $e_0 = |E(L_0, H_0)|$.
\end{defn}

\begin{claim} \label{G_0 charge}
$e_0 \leq 2(\ell+h)$.
\end{claim}
{\bf Proof.}
This is trivial if $h+\ell\leq 2$.
By definition, $L_0$ is  independent.
The claim follows by applying Lemma~\ref{co1}(i) for $A = L$ and $B=H$ for $h+\ell\geq 3$.
\qed

Let every vertex  $v\in V(J)$ have initial charge $d(v)$.
Our discharging has two rules:

{\bf Rule R1:} If a  vertex in $H_1$ has neighbors of degree $k-1$, then it keeps for itself charge $k-2/(k-1)$ and distributes the rest equally
among its neighbors of degree $k-1$. Otherwise, it keeps all its charge for itself.

{\bf Rule R2:} If a $K_{k-1}$-subgraph $X$ contains exactly $s\;$ $(k-1)$-vertices adjacent to a $(k-1)$-vertex $x$
outside of $X$ and not in a $K_{k-1}$, then each of these $s$ vertices gives charge $\frac{k-3}{s(k-1)}$ to $x$.

\begin{claim}\label{cl56}
Each vertex in $H_1$ gives to each neighbor of degree $k-1$ charge at least $\frac{1}{k-1}$.
\end{claim}
{\bf Proof.}
If $v \in H_1$, then $v$ gives  to each neighbor charge at least $\psi(d(v)):=\frac{d(v) - k + 2/(k-1)}{d(v)}$.
Since $\psi(x)$ is monotonically increasing for $x\geq k$,  $\psi(d(v))$ is minimized when $d(v) = k+1$.
Then each neighbor of $v$ of degree $k-1$ gets charge at least $(1 + 2/(k-1))/(k+1)=1/(k-1)$.
\qed

\begin{claim}\label{cl57}
Each vertex in $L_1$ has charge at least $k-2/(k-1)$.
\end{claim}
{\bf Proof.}
Let $v \in L_1$.
By Fact \ref{clusters partition} every vertex in $L \supseteq L_1$ is in exactly one cluster. 
Let $v$ be in a cluster $C$ of size $t$.  

By definition, if $v$ gives anything out to some vertex $x$ by Rule R2, then $v$ is in a $(k-1)$-clique $X$.
In this case, since $v$ has $k-2$ neighbors in $X$,  $x$ is the unique neighbor of $v$ outside of $X$. Also, by Claim \ref{unique cluster},  $C \subset X$
and by the definition of a cluster, each $w\in C$
is adjacent to $x$. So
\begin{equation}\label{j12}
\mbox{\em if $v$ gives anything out  by Rule 2, then it gives it at most once and for $s=t$.}
\end{equation}



{\bf Case 1:} $v$ is in a $(k-1)$-clique $X$ and $t \geq 2$.
Since  $v\notin H_1$, it can give charge only by Rule R2. Since $v$ is in a $(k-1)$-clique,
it may receive charge only by Rule R1.

Recall that $C \subset X$.
By Claim \ref{big neighbors (b)}, each vertex in $X-C$ has degree at least $k - 1 + t$.
Since $t \geq 2$, this means  $X-C \subseteq H_1$.
Furthermore, each vertex in $X-C$ has at least $k-2-t$ neighbors with degree at least $k$ (the other vertices of $X-C$).
Therefore each vertex $u \in X-C$ gives by Rule R1 charge at least $\frac{d(u) - k + 2/(k-1)}{d(u)-k+2+t}$  
to each neighbor of degree $k-1$.
Note that this function increases as $d(u)$ increases, so the charge is minimized when $d(u) = k-1+t$.
It follows that $u$ gives to $v$  charge at least $\frac{t-1+2/(k-1)}{2t+1}$.

So, $v$ receives from the vertices in $X-C$ by Rule R1 at least $(k-1-t)(\frac{t-1 + 2/(k-1)}{2t+1})$ and in view of~\eqref{j12}  gives
out by Rule R2 to at most one vertex and with $s=t$.
Thus $v$
 has charge at least $k-1 + (k-1-t)(\frac{t-1 + 2/(k-1)}{2t+1}) - \frac{k-3}{t(k-1)}$, which we claim is at least $k-2/(k-1)$.
Let
$$ g_1(t) =(k-1-t)((t-1)(k-1)+2) - (2t+1)(k-3)(1+\frac{1}{t}).
$$
We claim that $g_1(t) \geq 0$, which is equivalent to $v$ having charge at least $k - 2/(k-1)$.
Let
$$ \widetilde{g}_1(t) =(k-1-t)((t-1)(k-1)+2) - (2t+1)(k-3)(3/2).
$$
Note that $\widetilde{g}_1(t) \leq g_1(t)$ when $t \geq 2$, so it is enough to show that $\widetilde{g}_1(t) \geq 0$ on the appropriate domain.
Function $\widetilde{g}_1(t)$ is quadratic with a negative coefficient at $t^2$, so it suffices to check its values at the boundaries.
They are
$$ \widetilde{g}_1(2)  = (k-3)(k-6.5)
$$
and
\begin{eqnarray*}
  4 \widetilde{g}_1(\frac{k-1}{2}) & = & (k-1)\left( (k-3)(k-1) + 4 \right) - 6 k (k-3) \\
  		& = & k^3 - 11 k^2 + 29 k - 7 \\
  		& = &  (k-7)(k^2 - 4k + 1).
\end{eqnarray*}
Each of these values is non-negative when $k \geq 7$, and if $k=6$, then
 the case does not apply by Claim \ref{k=6}.

{\bf Case 2:} $t \geq 2$ and $v$ is not in a $(k-1)$-clique. Then $v$ cannot give charge by Rule R2
 and may receive charge by Rule R1. By Claim \ref{big neighbors (a)}, each neighbor of $v$ outside of $C$ has degree at least $k - 1 + t \geq k+1$ and is in $H_1$.
Therefore $v$ has charge at least $k-1 + (k-t)(\frac{t-1 + 2/(k-1)}{k-1+t})$.
We define
\begin{eqnarray*}
	g_2(t) &=& (k-t)(t-1+\frac{2}{k-1}) - \frac{k-3}{k-1}(k-1+t)\\
		&=& t(k-t) - 2(1 - \frac{2}{k-1})(k-1)\\
		& = & t(k-t) - 2(k - 3).
\end{eqnarray*}
Note that $g_2(t) \geq 0$ is equivalent to $v$ having charge at least $k - 2/(k-1)$.
The function $g_2(t)$ is quadratic with a negative coefficient at $t^2$, so it suffices to check
its values at the boundaries.
They are
$$ g_2(2) = 2(k-2) - 2(k-3) = 2
$$
and
$$ g_2(k-3) = (k-3)(3) - 2(k-3) = k-3.
$$
Each of these values is positive.

{\bf Case 3:} $t = 1$.
By the definition of $L_1$, $v$ is adjacent to at least one vertex $w$ with degree $k-1$.
Since $|C| = t = 1$ and so $C = \{v\}$,  $w \notin C$. So by Fact \ref{clusters partition},
$v$ and $w$ are in  different clusters. This means
 $N[w] \neq N[v]$.
If $v$ is not in a $(k-1)$-clique $X$, then by Claim \ref{adjacent k-1}, $w$ is in a $(k-1)$-clique and a cluster of size at least $2$.
In this case $v$ will receive charge $(k-3)/(k-1)$ in total from the cluster containing $w$ using Rule R2, and will not give away any charge.
Therefore we may assume that $v$ is in a $(k-1)$-clique $X$.

By Claim \ref{cliques have k+1}, there exists a $Y \subset X$ such that $|Y| \geq \frac{k-1}2$ and every vertex in $Y$ has degree at least $k+1$.
By Claim \ref{unique cluster},  none of the  vertices in $X - C = X - \{v\}$ is  in a cluster, which  by Fact \ref{clusters partition} 
means that every vertex in $X - \{v\}$ has degree at least $k$.
So each vertex in $Y$ has at least $k-3$ neighbors with degree at least $k$ (the vertices of $X$ besides $v$ and itself).
Therefore by Rule R1 each vertex $u \in Y$ donates charge at least $g_3(d(u))=\frac{d(u) - k + 2/(k-1)}{d(u)-k+3}$
 to each neighbor of degree $k-1$.
Since $d_3(d(u))$ increases as $d(u)$ increases,  the charge is minimized when $d(u) = k+1$.
It follows that $u$ gives to $v$ charge at least $d_3(k+1)=\frac{1+2/(k-1)}{4}$, and $v$ has charge at least
$$k-1 + \frac{k-1}2\left(\frac{1 + 2/(k-1)}{4}\right) = k + \frac{k-7}{8},
$$
which is at least $k-2/(k-1)$ when $k \geq 6$.
\qed

By the discharging rules and Claims~\ref{cl56} and~\ref{cl57},  after discharging,\\
a) the charge of each vertex in $H_1\cup L_1$ is at least $k-2/(k-1)$;\\
b) the charges of vertices in $H_0$ did not decrease;\\
c) along every edge from $H_1$ to $L_0$ the charge at least $1/(k-1)$ is sent.

Thus by Claim \ref{G_0 charge}, the total charge $F$ of the vertices in $H_0 \cup L_0$ is at least
$$kh+(k-1)\ell+\frac{1}{k-1}\left(\ell(k-1)-e_0\right) \geq k(h+\ell)-\frac{1}{k-1}2(h+\ell) = (h+\ell)\left(k-\frac{2}{k-1}\right),
$$
and so by a), the total charge of all the vertices of $J$ is at least $n\left(k-\frac{2}{k-1}\right)$.
Therefore the degree sum of $J$ is at least $n\left(k-\frac{2}{k-1}\right) = \left(\frac{(k+1)(k-2)}{k-1}\right)n$, i. e., $|E(J)| \geq \left(\frac{(k+1)(k-2)}{2(k-1)}\right)n$.
\qed

\section{Sharpness}

First we prove Corollary \ref{new tightness}, and then we will construct sparse $3$-connected $k$-critical graphs.
As it was pointed out in the introduction, Construction~\ref{3ConnConst} and
infinite series of $3$-connected sparse $4$-~and $5$-critical graphs are due to Toft~\cite{T12} (based on~\cite{Toft2}).

{\bf Proof of Corollary \ref{new tightness}.}
By  (\ref{upper f_k}), if we construct an $n_0$-vertex $k$-critical graph for which our lower bound on $f_k(n_0)$
is exact, then the bound on $f_k(n)$ is exact for every $n$ of the form $n_0+s(k-1)$.
So, by Corollary \ref{Ore Cor}, we only need to construct
\begin{itemize}
	\item a $5$-critical $7$-vertex graph with $\left\lceil 15 \frac12\right\rceil = 16$ edges,
	\item a $5$-critical $8$-vertex graph with $\left\lceil 17 \frac34\right\rceil = 18$ edges,
	\item a $6$-critical $10$-vertex graph with $\left\lceil 27 \frac15\right\rceil = 28$ edges,
	\item a $6$-critical $12$-vertex graph with $\left\lceil 32 \frac45\right\rceil = 33$ edges, and
	\item a $7$-critical $14$-vertex graph with $\left\lceil 45 \frac13\right\rceil = 46$ edges.
\end{itemize}

\begin{figure}[htbp]
\begin{center}
\includegraphics[height=4cm]{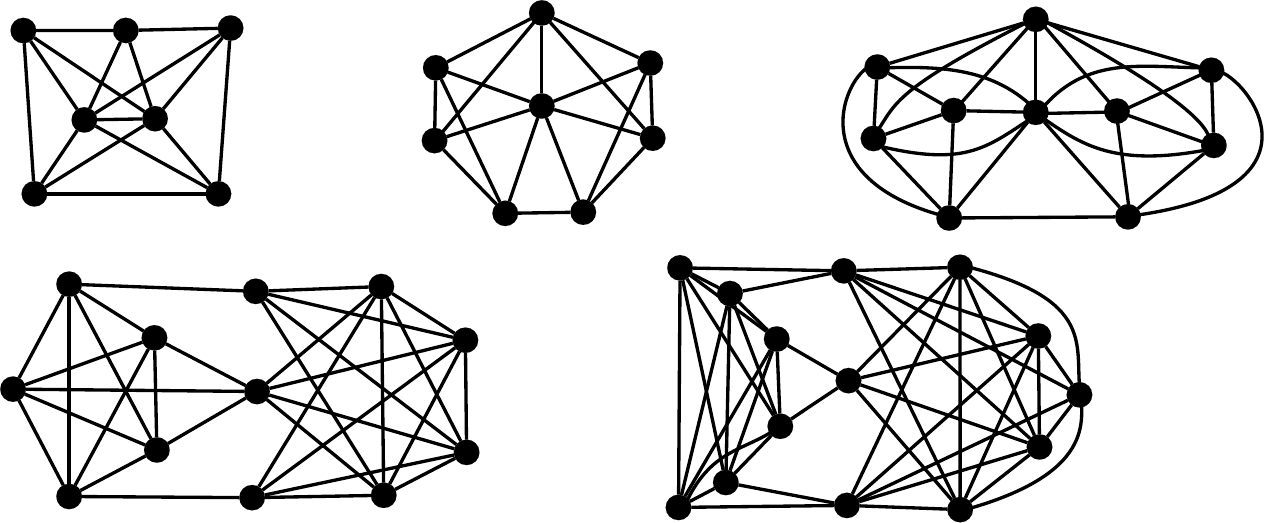}

\caption{Minimal $k$-critical graphs.}
\label{new tightness examples}
\end{center}
\end{figure}

These graphs are presented in Figure \ref{new tightness examples}.
\qed

\begin{const}[Toft~\cite{T12}] \label{3ConnConst}
Let $G$ be a $k$-critical graph, $e = uv \in E(G)$, and $w \in V(G)-\{u,v\}$ be such that for
all $(k-1)$-colorings $\phi$ of $G-e$, $\phi(w) = \phi(u) = \phi(v)$.
Let $S_1 \cup S_2 \cup S_3$ be a partition of the vertex set $X$ of a copy of $K_{k-1}$  such that each $S_i$ is non-empty.
We define  $G'$ as follows:  $V(G') = V(G) \cup V(X)$ and $E(G') = (E(G)-e) \cup E(X) \cup E'$, where
$$ E' = \{ua : a \in S_1\}  \cup \{vb : b \in S_2\} \cup \{wc : c \in S_3\}.$$
\end{const}

\begin{claim} \label{InductionStep}
If $k\geq4$, $G$ is a $3$-connected $k$-critical graph and $G'$ is created from $G$ 
via  Construction \ref{3ConnConst}, then $G'$ is a $3$-connected $k$-critical graph.
\end{claim}
{\bf Proof.}
We will use the  names and definitions from Construction \ref{3ConnConst}.

If there exists a $(k-1)$-coloring $\phi$ of $G'$, then all $k-1$ colors must appear on $X$.
Then $\phi(u)$ appears on a vertex in $S_2$ or $S_3$.
But then either $\phi(v)\neq \phi(u)$ or $\phi(w) \neq \phi(u)$, which contradicts the assumptions of Construction \ref{3ConnConst}.
So $\chi(G') \geq k$.

Suppose there exists an $f \in E(G')$ such that $\chi(G' - f) \geq k$.
If $f \in E(G)$, then let $\phi_1$ be a $(k-1)$-coloring of $G-f$.
Because $e \in E(G)-f$, $\phi_1(u) \neq \phi_1(v)$, and so $\phi_1$ extends easily to $G'-f$.
If $f \subset X$, then a $(k-1)$-coloring of $G-e$ can be extended to $G'-f$, because $X$ can
be colored with $k-2$ colors, while $N(X) = \{u,v,w\}$ is colored with $1$ color.
If $f \in E'$, then a $(k-1)$-coloring of $G-e$ extends to $G'-f$, because the unique color on $\{u,v,w\}$ can be given to $f \cap X$.
Therefore $G'$ is $k$-critical.

Suppose now that there exists a set $S$ such that $|S| < 3$ and there are nonempty $A$ and $B$ such that $E(A,B) = \emptyset$ and $A \cup B \cup S = V(G')$.
Because $k$-critical graphs are $2$-connected, $|S| = 2$.  By Fact~\ref{fa7}, $S$ is independent.
Because $X$ is a clique, without loss of generality $X \subsetneq A \cup S$.
By construction,  $A-X \neq \emptyset$, so $S$ also separates $G-e$.
Since $\kappa(G) \geq 3$,  $S\cap X=\emptyset$ and $e$ has an endpoint in each component of $G-S-e$.
So we may assume that $u\in A$ and $v\in B$. Since $v$ has a neighbor in $S_2\subset X$, this contradicts to
$S\cap X=\emptyset$ and  $X  \subsetneq A \cup S$.
\qed

The assumptions in Construction \ref{3ConnConst} are strong.
Most edges $e$ in $k$-critical graphs do not have such a vertex $w$, and some $k$-critical graphs do not have any edge-vertex pairs $(e,w)$ that satisfy the assumptions.
We will construct an infinite family $\mathbb{G}_k$ of sparse $3$-connected graphs that do satisfy the assumptions.

The family is generated for each $k$ by finding a small $3$-connected $k$-critical graph $G'_k$ such that $\rho_k(G'_k) = y_k$.
We will describe a subgraph $H_k' \leq G_k'$ with two vertices, $u$ and $w$, such that in any $(k-1)$-coloring $\phi'$ of $H_k'$, $\phi'(u) = \phi'(w)$.
Construction \ref{3ConnConst} can then be applied to $G_k'$, using any edge $e$ incident to $u$ that is not in $H_k'$ and not incident to $w$.
Because Construction \ref{3ConnConst} does not decrease the degree of $u$, this process can be iterated indefinitely to populate $\mathbb{G}_k$.

Note that Construction \ref{3ConnConst} adds the same number of vertices and edges as \Ore-composition with $G_2 = K_k$.
Therefore every graph $G \in \mathbb{G}_k$ has $\rho_k(G) = y_k$.
Furthermore, $G$ is also $k$-critical and $3$-connected, and therefore not $k$-Ore.
This implies the sharpness of Theorem \ref{ext}.

\begin{figure}[htbp]
\begin{center}
\includegraphics[height=4cm]{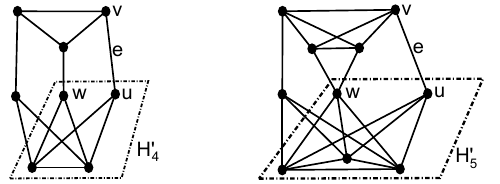}

\caption{Graphs $G'_4$ and $G'_5$, with  substructures
 labeled  for constructing $\mathbb{G}_4$ and $\mathbb{G}_5$.}
\label{small 3conn graphs}
\end{center}
\end{figure}

All that is left is to find suitable graphs for $G'_k$ and $H'_k$.
Figure \ref{small 3conn graphs} illustrates $G'_4$ and $G'_5$.
We will need a second construction for larger $k$.

\begin{const}\label{con2}
Fix a $t$ such that $1 \leq t < k/2$.
Let
$$V(H_{k,t}) = \{u_1, u_2, \ldots, u_{k-1}, v_1, v_2, \ldots, v_{k-1}, w\}$$
 and
$$ E(H_{k,t}) = \{u_iu_j:1 \leq i < j \leq k-1\} \cup \{v_iv_j:1 \leq i < j \leq k-1\} \cup \{u_iv_j: i,j \leq t\} $$
$$ \cup \{w u_i: i > t\} \cup \{w v_i: i > t\}.$$
\end{const}

By construction, $H_{k,1}$ is a $k$-Ore graph, $H_{k,t}$ is $k$-critical, $\kappa( H_{k,t} ) = t+1$, $|V(H_{k,t})| = 2k-1$, and $|E(H_{k,t})| = k(k-1) - 2t + t^2$.
Moreover, $\rho_k(H_{k,2}) = y_k$.
For $k \geq 6$, we choose $G'_k = H_{k,2}$.
We will next find $H'_k$ for $k \geq 6$, which will complete the argument.

\begin{claim}
Let $H'_k = H_{k,2} - u_1v_1$.
Then in every $(k-1)$-coloring $\phi'$ of $H'_k$, $\phi'(u_1) = \phi'(v_1)=\phi'(w)$.
\end{claim}
{\bf Proof.}
Let $\phi'$ be a $(k-1)$-coloring of $H'_k$.
Note that each of the $(k-1)$ colors appears both, on $\{u_1, u_2, \ldots , u_{k-1}\}$ and 
on $\{v_1, v_2, \ldots , v_{k-1}\}$.
Then $\phi'(w)$ appears on a vertex $a \in \{u_1, u_2\}$ and again on a vertex $b \in \{v_1, v_2\}$.
So $ab \notin E(G)$, which implies that $a = u_1$ and $b=v_1$.
\qed

\begin{figure}[htbp]
\begin{center}
\includegraphics[height=4cm]{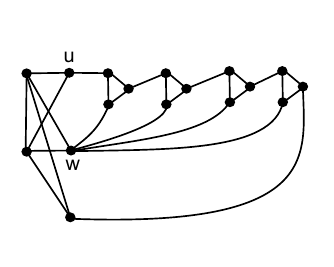}

\caption{An example of a graph in $\mathbb{G}_4$.}
\label{G44}
\end{center}
\end{figure}

\section{Algorithm}
The proof of Theorem \ref{k-critical} was constructive, and provided an algorithm for $(k-1)$-coloring of sparse graphs.
Let $P_k(G)$ be the minimum of $\rho_{k,G}(W)$ over all $W\subseteq V(G)$ with $2\leq |W|$.

\vspace{.1in}


\begin{theorem}[\cite{KY}] \label{algo}
If $k \geq 4$, then every $n$-vertex graph $G$ with $P_k(G) > k(k-3)$ can be $(k-1)$-colored in $O(k^{3.5}n^{6.5}\log(n))$ time.
\end{theorem}

We present below a polynomial-time algorithm for checking whether a given graph is a $k$-Ore graph.
Together with an analog of the algorithm in Theorem~\ref{algo} that uses
the proof of
Theorem~\ref{ext} instead of Theorem~\ref{k-critical},
it would yield a polynomial-time algorithm that for every $n$-vertex graph $G$ with $P_k(G) > y_k$ either finds
a $(k-1)$-coloring of $G$ or finds a subgraph of $G$ that is a $k$-Ore graph.

Our algorithm to determine  whether an $n$-vertex graph $G$ is $k$-Ore is simple:\\

0. If $G$ is $K_k$, return ``yes.''\\

1. Check whether $n\equiv 1 \,(\mod k-1)$ and $|E(G)|=\frac{(k+1)(k-2)|V(G)|-k(k-3)}{2(k-1)}$.
If not, then return ``no.''\\

2. Check whether the connectivity of $G$ is exactly $2$. If not, then return ``no.''
Otherwise, choose a separating set $\{x,y\}$.\\

3. If $G-x-y$ has more than two components or $xy\in E(G)$,  then return ``no.''
Otherwise, let $A$ and $B$ be the vertex sets of the two components of $G-x-y$.
If $\{|A| \,(\mod k-1),|B| \,(\mod k-1)\}\neq \{k-2,0\}$, then return ``no''. Otherwise,
rename $A$ and $B$ so that $|A| \,(\mod k-1)= k-2$ and $|B| \,(\mod k-1)=0$.\\

4.  Create graphs $\widetilde G(x,y)$ and $\check G(x,y)$ as defined in Fact \ref{f2}.
  Recurse on each of $\widetilde G(x,y)$ and $\check G(x,y)$.
  If at least one recursion call returns ``no,'' then return ``no.''  Otherwise, return ``yes.''\\

The longest procedure in this algorithm is checking whether the connectivity of $G$ is exactly $2$
at Step 2, which has complexity $O(kn^3)$ because $|E(G)| \leq kn/2$. 
And it will be called fewer than $2n/(k-2)$ times. 
So the overall complexity is at most $O(n^4)$.

\paragraph{Acknowledgment.} We thank Michael Stiebitz and Bjarne Toft for helpful discussions.
We also thank the referees for many helpful comments  significantly improving the presentation of the paper.

\end{document}